\documentclass[11pt]{article}
\usepackage[dvips]{graphicx}
\usepackage{latexsym}
\usepackage{amssymb}

\newcommand{\C}{\mathbb{C}}
\newcommand{\CP}{\mathbb{CP}}

\newcommand{\Z}{\mathbb{Z}}
\newcommand{\X}{\mathbb{X}}

\renewcommand{\d}{\mathrm{d}}

\newcommand{\koniec}{\begin{flushright}  $\Box $ \end{flushright}}
\def\be{\begin{equation}}
\def\ee{\end{equation}}

\def\Sm{\Sigma}

\def\Om{\Omega}
\def\Th{\Theta}
\def\G{\Gamma}
\def\O{\cal O}
\def\om{\omega}
\def\Do{\p _0}
\def\dd{\p _2}

\def\wt{\overline{w}}
\def\zt{\overline{z}}

\def\tw{\tilde w}
\def\tz{\tilde z}

\def\p{\partial}
\def\pt{\tilde \p}
\def\do{\p _0}
\def\dd{\p _2}

\newcommand{\hook}{{\setlength{\unitlength}{11pt}   % adjust pt size here
                   \begin{picture}(.833,.8)
                   \put(.15,.08){\line(1,0){.35}}
                   \put(.5,.08){\line(0,1){.5}}
                   \end{picture}}}
\def\a{\alpha}
\def\l{\lambda}
\def\Pt{{\cal PT}}
\def\O{{\cal O}}
\def\at{\tilde \a}
\def\dom{\delta \Om}
\def\tom{\tilde{\omega}}
\def\ddot{\delta \Theta}

\def\TH{
\left (
\begin{array}{cc}
\p _y& \p _w+\Th_{yy}\p _x-\Th _{xy}\p_y  \\
-\p_x&\p _z-\Th _{xy}\p _x + \Th _{xx}\p_y
\end{array}
\right ) 
}
\def\OM{
\left (
\begin{array}{cc}
\Om _{w\wt}\p _{\zt}-\Om _{w\zt}\p _{\wt}&\p_w \\
\Om _{z\wt}\p _{\zt}-\Om _{z\zt}\p _{\wt}&\p_z
\end{array}
\right ) 
}

\newtheorem{theo}{Theorem}[section] 
\newtheorem{prop}[theo]{Proposition}  
\newtheorem{lemma}[theo]{Lemma}
\newtheorem{defi}[theo]{Definition}

\begin{document}

\title{Hyper-K\"ahler Hierarchies and their twistor theory}
\author{Maciej  Dunajski\thanks{email: dunajski@maths.ox.ac.uk},
Lionel J. Mason\\ The Mathematical Institute,
24-29 St Giles, Oxford OX1 3LB, UK}  
\date{}
\maketitle
\noindent
\abstract{A twistor construction of the hierarchy associated with the
  hyper-K\"ahler equations on a metric (the anti-self-dual Einstein
  vacuum equations, ASDVE, in four dimensions) is given. The recursion
  operator $R$ is constructed and used to build an
  infinite-dimensional symmetry algebra and in particular higher flows
  for the hyper-K\"ahler equations.  It is shown that $R$ acts on the
  twistor data by multiplication with a rational function.  The
  structures are illustrated by the example of the Sparling-Tod
  (Eguchi-Hansen) solution.

  An extended space-time ${\cal N}$ is constructed whose extra
  dimensions correspond to higher flows of the hierarchy.  It is shown
  that ${\cal N}$ is a moduli space of rational curves with normal
  bundle ${\cal O}(n)\oplus{\cal O}(n)$ in twistor space and is
  canonically equipped with a Lax distribution for ASDVE hierarchies.
  The space ${\cal N}$ is shown to be foliated by four dimensional
  hyper-K{\"a}hler slices.
  
  The Lagrangian, Hamiltonian and bi-Hamiltonian formulations of the
  ASDVE in the form of the heavenly equations are given.  The
  symplectic form on the moduli space of solutions to heavenly
  equations is derived, and is shown to be compatible with the
  recursion operator. }

\section{Introduction}
Roger Penrose's twistor theory gives rise to correspondences between
solutions to differential equations on the one hand and unconstrained
holomorphic geometry on the other.  The two most prominent systems of
nonlinear equations which admit such correspondences are the
anti-self-dual vacuum Einstein equations (ASDVE) \cite{Pe76} which 
in Euclidean
signature determine hyper-K\"ahler metrics, and the
anti-self-dual Yang--Mills equations (ASDYM) \cite{Wa77}.  Richard
Ward \cite{Wa85} observed that many lower-dimensional integrable
systems are symmetry reductions of ASDYM. This has led to an
overview of the theory of integrable systems \cite{MW95}, which
provides a classification of those lower-dimensional integrable
systems that arise as reductions of the ASDYM equations and a
unification of the theory of such integrable equations as symmetry
reduced versions of the corresponding theory of the ASDYM equations.
In \cite{MW95}, Lagrangian and Hamiltonian frameworks for ASDYM were
described together with a recursion operator.  This leads to the
corresponding structures for symmetry reductions of the ASDYM
equations.  

In this paper we investigate these structures for the second important
system of equations---the ASDVE or hyper-K\"ahler equation (this system
also admits known integrable systems as symmetry reductions
\cite{DMW95}).  We shall give a twistor-geometric construction of the
hierarchies associated to the ASDVE in the `heavenly' forms due to
Pleba\'nski \cite{Pl75}.  In this context it is more natural to work
with complex (holomorphic) metrics on complexified space-times and so
we use the term ASDVE equations rather than hyper-K\"ahler equations.
Our considerations will generally be local in space-time which will be
understood to be a region in $\C^4$.

In Section \ref{preliminaries} we summarise the twistor
correspondences for flat and curved spaces.  We establish a spinor
notation (which will not be essential for the subsequent sections) and
recall basic facts about the ASD conformal condition and the geometry
of the spin bundle.  In Section \ref{Heavenly_Hierarchies} the
recursion operator $R$ for the ASDVE is constructed as an integro-
differential operator mapping solutions to the linearised heavenly
equations to other solutions.  We then use this to give an alternate
development of the twistor correspondence by using $R$ to build a
family of foliations by twistor surfaces.  We show that $R$
corresponds to multiplication of the twistor data by a given twistor
function.  We then analyse the hidden symmetry algebra of the ASDVE,
and use the recursion operator to construct Killing spinors.  We
illustrate the ideas using the example of the Sparling--Tod solution
and show how $R$ can be used to construct rational curves with normal
bundle ${\cal O}(1)\oplus{\cal O}(1)$ in the associated twistor space.

In Section \ref{Heavenly_h} we give the twistor construction for the
ASDVE hierarchies.  The higher commuting flows can be thought of as
coordinates on an extended space-time.  This extended space-time has a
twistor correspondence: it is the moduli space of rational curves with
normal bundle ${\cal O}(n)\oplus{\cal O}(n)$ in a twistor space.  This
moduli space is canonically equipped with the Lax distribution for
ASDVE hierarchies, and conversely that truncated hierarchies admit a
Lax distribution that gives rise to such a twistor space.  The Lax
distribution can be interpreted as a connecting map in a long exact
sequence of sheaves.  In Section \ref{hformalism} we investigate
the Lagrangian and Hamiltonian formulations of heavenly equations.
The symplectic form on the moduli space of solutions to heavenly
equations will be derived, and is shown to be compatible with the
recursion operator.

We end this introduction with some bibliographical remarks.
Significant progress towards understanding the symmetry structure of
the heavenly equations was achieved by Boyer and Pleba\'nski
\cite{BP77,BP85} who obtained an infinite number of conservation laws
for the ASDVE equations and established some connections with the
nonlinear graviton construction.  Their results were later extended in
papers of Strachan \cite{St95} and Takasaki \cite{Ta89,Ta90}.  The present
work is an extended version of \cite{DM96,DM97,DM97c}.

\section{Preliminaries}
\label{preliminaries}

\subsection{Spinor notation}
We work in the holomorphic category with complexified space-times:
thus space-time ${\cal M}$ is a complex four-manifold equipped with a
holomorphic metric $g$ and compatible volume form $\nu$.

In four complex dimensions orthogonal transformations decompose into
products of ASD and SD rotations 
\be
\label{basicisom}
SO(4, \C)=(SL(2, \C)\times \widetilde{SL}(2, \C))/\Z_2.  \ee The
spinor calculus in four dimensions is based on this isomorphism.  We
use the conventions of Penrose and Rindler \cite{PR86}.  Indices will
generally be assumed to be concrete unless stated otherwise:
%if we work in any of the heavenly frames and otherwise abstract: 
$a,b,\ldots$, $a=0,1\ldots 3$ are four-dimensional space-time indices
and $A, B, \ldots, A', B', \ldots$, $A=0,1$ etc.\ are two-dimensional
spinor indices.  The tangent space at each point of ${\cal M}$ is
isomorphic to a tensor product of the two spin spaces  
\be
\label{isomorphism}
T^a{\cal M}=S^A\otimes S^{A'}.
\ee
The complex Lorentz transformation $V^a\longrightarrow
{\Lambda^a}_bV^b$, ${\Lambda^a}_b{\Lambda^c}_d g_{ac}=g_{bd}$,
is equivalent to the composition of the SD and the ASD rotation
\[
V^{AA'}\longrightarrow {\lambda^A}_BV^{BB'}{\lambda^{A'}}_{B'},
\]
where ${\lambda^A}_B$ and ${\lambda^{A'}}_{B'}$ are elements of $SL(2, \C)$
and $\widetilde{SL}(2, \C)$.

Spin dyads $(o^A, \iota^{A})$ and $(o^{A'}, \iota^{A'})$ span $S^A$
and $S^{A'}$ respectively.  The spin spaces $S^A$ and $S^{A'}$ are
equipped with symplectic forms $\varepsilon_{AB}$ and
$\varepsilon_{A'B'}$ such that
$\varepsilon_{01}=\varepsilon_{0'1'}=1$.  These anti-symmetric objects
are used to raise and lower the spinor indices. We shall use 
normalised spin frames so  that
\[
o^B\iota^C-\iota^Bo^C=\varepsilon^{BC},\;\;\;
o^{B'}\iota^{C'}-\iota^{B'}o^{C'}=\varepsilon^{B'C'}.
\] 
 Let $e^{AA'}$ be a null tetrad of 1-forms on ${\cal M}$ and let
$\nabla_{AA'}$ be the frame of dual vector fields. 
The orientation is given by fixing the volume form
\[
\nu=e^{01'}\wedge e^{10'}\wedge e^{11'}\wedge e^{00'}.
\]
Apart from orientability,
${\cal M}$ must satisfy some other topological restrictions for the
global spinor fields to exist. We shall not take them into
account as we work locally in ${\cal M}$.

The local basis $\Sm^{AB}$ and $\Sm^{A'B'}$ of  spaces of ASD and SD
two-forms are defined by
\be\label{sdforms}
e^{AA'}\wedge e^{BB'}=\varepsilon^{AB}\Sm^{A'B'}+\varepsilon^{A'B'}\Sm^{AB}. 
\ee
The first Cartan structure equations are 
\[
\d e^{AA'}=e^{BA'}\wedge{\Gamma^{A}}_{B}+e^{AB'}\wedge{\Gamma^{A'}}_{B'},
\]
where $\Gamma_{AB}$ and $\Gamma_{A'B'}$ are the  $SL(2, \C)$
and $\widetilde{SL}(2, \C)$
spin connection one-forms.  They are symmetric in their indices, and
\[
 \Gamma_{AB}=
\Gamma_{CC'AB}e^{CC'},\;\;\Gamma_{A'B'}=\Gamma_{CC'A'B'}e^{CC'}
,\;\;\; \Gamma_{CC'A'B'}=o_{A'}\nabla_{CC'}\iota_{B'}-
\iota_{A'}\nabla_{CC'}o_{B'}.
\]
The curvature of the spin connection
\[
{R^A}_B=\d{\Gamma^A}_B+{\Gamma^A}_C\wedge{\Gamma^C}_B
\]
decomposes as
\[
{R^A}_B={C^A}_{BCD}\Sm^{CD}+(1/12)R{\Sm^{A}}_{B}+{\Phi^A}_{BC'D'}\Sm^{C'D'},
\]
and similarly  for ${R^{A'}}_{B'}$. Here $R$ is the Ricci scalar, 
$\Phi_{ABA'B'}$ is the trace-free part of the Ricci tensor $R_{ab}$,
and $C_{ABCD}$ is the ASD part of the Weyl tensor
\[
C_{abcd}=\varepsilon_{A'B'}\varepsilon_{C'D'}C_{ABCD}+
\varepsilon_{AB}\varepsilon_{CD}C_{A'B'C'D'}.
\]
%%%%%%%%%%%%%%%%%%%%%%%%%%%%%%%%%%%%%%%%%%%%%%%%%%%%%%%%%%%%%%%%%%%%%%%%%%%%%%%

\subsection{The flat twistor correspondence}
The flat twistor correspondence is a correspondence between points in
complexified Minkowski space, $\C^4$ (or its conformal
compactification) and holomorphic lines in $\CP^3$.  

The flat twistor correspondence has an invariant formulation in terms
of spinors.  A point in $\C^4$ has position vector with coordinates
$(w, z, x, y)$.  The isomorphism (\ref{isomorphism}) is realised by
\[
x^{AA'}:=
\left (
\begin{array}{cc}
y&w\\
-x&z
\end{array}
\right ),\;\;\;\;\mbox{so that}\;\; 
g=\varepsilon_{AB}\varepsilon_{A'B'}\d x^{AA'}\d x^{BB'}. 
\]
A two-plane in  $\C^4$ is null if $g(X, Y)=0$ for every pair $(X, Y)$
of vectors tangent to it. The null planes can be self-dual (SD) or 
anti self-dual (ASD), depending on whether the tangent bi-vector
$X\wedge Y$ is SD or ASD. The SD null planes are called $\a$-planes.
The $\a$-planes passing through a point in $\C^4$ are parametrised by
$\l=\pi_{0'}/\pi_{1'}\in\CP^1$. Tangents to $\a$-planes are spanned by
two vectors  
\be
\label{flatlax}
L_A=\pi^{A'}\frac{\p}{\p x^{AA'}}
\ee 
which form the kernel of $\pi_{A'}\pi_{B'}\Sigma^{A'B'}$
The set of all $\a$-planes is called a projective twistor space and
denoted ${\cal PT}$. For $\C^4$ it is a three-dimensional complex manifold
bi-holomorphic to $\CP^3-\CP^1$.

The five complex dimensional correspondence space ${\cal
F}:=\C^4\times\CP^1$ fibres over $\C^4$ by $(x^{AA'}, \l) \rightarrow
x^{AA'}$ and over ${\cal PT}$ with fibres spanned by $L_A$. Twistor
functions (functions on ${\cal PT}$) pull back to functions on ${\cal
F}$ which are constant on $\a$-planes, or equivalently satisfy
$L_Af=0$.

Twistor space can be covered by two coordinate patches $U$
and $\widetilde U$, where $U$ is a complement of $\l=\infty$ and
$\widetilde U$ is a compliment of $\l=0$. If $(\mu^0, \mu^1, \l)$ are
coordinates on $U$ and $(\tilde{\mu}^0, \tilde{\mu}^1, \tilde{\l})$
are coordinates on $\widetilde{U}$ then on the overlap
\[
\tilde{\mu}^0={\mu}^0/\l,\;\;\tilde{\mu}^1={\mu}^1/\l,\;\;\tilde{\l}=1/\l.
\]
The local coordinates $(\mu^0, \mu^1, \l)$ 
on ${\cal PT}$ pulled back to ${\cal F}$
are
\be
\label{tcoord}
\mu^0=w+\l y,\;\;\;\mu^1= z-\l x,\;\;\; \l.
\ee
We can introduce homogeneous coordinates on the twistor space
\[
(\om^A, \pi_{A'})=(\om^0, \om^1, \pi_{0'}, \pi_{1'}):= (\mu^0\pi_{1'}
, \mu^1\pi_{1'} ,\l\pi_{1'} ,\pi_{1'} ).
\]

The point $x^{AA'}\in \C^4$ lies on the $\alpha$-plane corresponding
to the twistor $(\om^A,\pi_{A'})\in{\cal PT}$ iff 
\be
\label{basiceq}
\om^A=x^{AA'}\pi_{A'}.  
\ee 
For $\pi_{A'}\neq 0$ and $(\om^A, \pi_{A'})$ fixed, The solution to
(\ref{basiceq}) is a complex two plane with tangent vectors of the form
$\pi^{A'}\nu^A$ for all $\nu^A$.  Alternatively, if we fix $x^{AA'}$, then
(\ref{basiceq}) defines a rational curve, $\CP^1$, in ${\cal PT}$ with
normal bundle ${\cal O}(1)\oplus {\cal O}(1)$.\footnote{ Here ${\cal
O}(n)$ denotes the line bundle over $\CP^1$ with transition functions
$\l^{-n}$ from the set $\l\neq\infty$ to $\l\neq 0$ (i.e.\ Chern class
$n$).}
Kodaira theory guarantees that the family of such rational curves in
${\cal PT}$ is four complex dimensional.  There is a canonical
(quadratic) conformal structure $\d s^2$ on $\C^4$: the points $p$ and
$q$ are null separated with respect to $\d s^2$ in $\C^4$ iff the
corresponding rational curves $l_p$ and $l_q$ intersect in ${\cal
PT}$ at one point.

%%%%%%%%%%%%%%%%%%%%%%%%%%%%%%%%%%%%%%%%%%%%%%%%%%%%%%%%%%%%%%%%%%%%%%%%%%%%

\subsection{Curved twistor spaces and the geometry of the primed spin
bundle.}

Given a complex four-dimensional manifold ${\cal M}$ with curved
metric $g$, a twistor in ${\cal M}$ is an $\a$-surface, i.e.\ a null 
two-dimensional surface whose tangent space at each point is an $\a$
plane.  There are Frobenius integrability conditions for the existence
of such $\a$-surfaces through each $\a$-plane element at each point
and these are equivalent, after some calculation, to the vanishing of
the self-dual part of the Weyl curvature, $C_{A'B'C'D'}$.  Thus, given
$C_{A'B'C'D'} =0$, we can define a twistor space ${\cal PT}$ to be the
three complex dimensional manifold of $\a$-surfaces in ${\cal M}$.  If
$g$ is also Ricci flat then ${\cal PT}$ has further structures which
are listed in the Nonlinear Graviton Theorem:

\begin{theo}[Penrose \cite{Pe76}]
\label{Penrose}
There is a 1-1 correspondence between complex ASD vacuum metrics
on complex four-manifolds 
and three dimensional complex manifolds ${\cal PT}$ such that
\begin{itemize}
\item There exists a holomorphic projection $\mu:{\cal
PT}\longrightarrow \CP^1$
\item ${\cal PT}$ is equipped with a four complex parameter
family of sections of $\mu$ each with a normal bundle $
{\cal O}(1)\oplus{\cal O}(1)$, (this will follow from the existence of
one such curve by Kodaira theory),
\item Each fibre of $\mu$ has a symplectic  structure
$
\Sm_{\l} \in \Gamma(\Lambda^2(\mu^{-1}(\lambda))\otimes {\cal O}(2)),
$
where $\l\in \CP^1$.
\end{itemize}
\end{theo}
To obtain real metrics on a real 4-manifold, we can require further that the
twistor space admit an anti-holomorphic involution.

The correspondence space ${\cal F}={\cal M}\times\CP^1$ is
coordinatized by $(x,\lambda)$, where $x$ denotes the coordinates on
$\cal M$ and $\l$ is the coordinate on $\CP^1$ that parametrises the
$\a$-surfaces through $x$ in $\cal M$.  We represent $\cal F$ as the
quotient of the primed-spin bundle $S^{A'}$ with fibre coordinates
$\pi_{A'}$ by the Euler vector field $\Upsilon=\pi^{A'}/\p
\pi^{A'}$. We relate the fibre coordinates to $\lambda$ by
$\lambda=\pi_{0'}/\pi_{1'}$.
A form with values in the line bundle ${\cal O}(n)$ on $\cal F$ can be
represented by a homogeneous form $\alpha$ on the non-projective spin
bundle satisfying
\[
\Upsilon\hook\alpha=0\, , \qquad {\cal L}_{\Upsilon} \alpha=n\alpha .
\]

The space $\cal F$ possesses a natural two dimensional distribution
called the twistor distribution, or Lax pair, to emphasise the analogy
with integrable systems.  The Lax pair on ${\cal F}$ arises as the
image under the projection $TS^{A'}\longrightarrow T{\cal F}$ of the
distribution spanned by
\[
\pi^{A'}\p_{AA'}+\Gamma_{AA'B'C'}\pi^{A'}\pi^{B'}\frac{\p}{\p
\pi_{C'}}
\]
on $TS^{A'}$ where the $\p_{AA'}$ are a null tetrad for the metric on
$\cal M$, and $\Gamma_{AA'B'C'}$ are the components of the spin
connection in the associated spin frame
($\p_{AA'}+\Gamma_{AA'B'C'}\pi^{B'}\frac{\p}{\p\pi_{C'}}$ is the
horizontal distribution on $S_{A'}$).  We can also represent the Lax
pair on the projective spin bundle by 
\footnote{ Various powers of $\pi_{1'}$ in formulae like
(\ref{laxpair}) guarantee the correct homogeneity. We usually shall
omit them when working on the projective spin bundle.  In a projection
$S^{A'}\longrightarrow {\cal F}$ we shall use the replacement formula
\be \label{procedura} 
\frac{\p}{\p
\pi^{A'}}\longrightarrow\frac{\pi_{A'}}{{{\pi_{1'}}^2}}\p_{\lambda}.
\ee This is because (on functions of $\lambda$) \[
\frac{\p}{\p\pi_{A'}}\Big(\frac{\pi_{0'}}{\pi_{1'}}\Big)=
\frac{\pi_{1'}o^{A'}-\pi_{0'}\iota^{A'}}{{\pi_{1'}}^2}
=\frac{\pi^{A'}}{{\pi_{1'}}^2}.  \] } \be
\label{laxpair}
L_A=(\pi_{1'}^{-1})(\pi^{A'}\p_{AA'}+f_A\p_{\lambda}),\;\;\;\;\mbox{
  where }\;\;
f_A=(\pi_{1'}^{-2})\Gamma_{AA'B'C'}\pi^{A'}\pi^{B'}\pi^{C'}.  
\ee 
The integrability of the twistor distribution is equivalent to
$C_{A'B'C'D'}=0$, the vanishing of the self-dual Weyl spinor. When the
Ricci tensor vanishes also, a covariant constant primed spin frame can
be found so that $\Gamma_{AA'B'C'}=0$.  We assume this from now on.

The projective twistor space $\cal PT$ arises as a quotient of $\cal
F$ by the twistor distribution.  With the Ricci flat condition, the
coordinate $\lambda$ descends to twistor space and $\pi_{A'}$ descends
to the non-projective twistor space.  It can be covered by two sets,
$U=\{|\lambda|< 1+\epsilon\}$ and
$\tilde{U}=\{|\lambda|>1-\epsilon\}$.  On the non-projective space we
can introduce extra coordinates $\om^A$ of homogeneity degree one so
that $(\omega^A, \pi_{A'}), \pi_{A'}\neq\iota_{A'}$ are homogeneous
coordinates on $U$ and similarly $(\tilde{\om}^A,\pi_{A'}), \pi_{A'}
\neq o_{A'})$ on $\tilde{U}$.  The twistor space $\cal PT$ is then
determined by the transition function ${\tilde{\om}}^B=
{\tilde{\om}}^B(\om^A, \pi_{A'})$ on $U\cap \tilde{U}$.

The correspondence space has the alternate definition
\[
{\cal F}={\cal PT}\times {\cal M}|_{Z\in l_x}= {\cal M}\times\CP^1
\] 
where $l_x$ is the line in $\cal PT$ that corresponds to $x\in {\cal
M}$ and $Z\in\cal PT$ lies on $l_x$.  This leads to a double fibration
\be
\label{doublefib}
{\cal M}\stackrel{p}\longleftarrow 
{\cal F}\stackrel{q}\longrightarrow {\cal PT}.
\ee

The existence of $L_A$ can also be deduced directly from the
correspondence.  From \cite{Pe76}, points in ${\cal M}$ correspond to
rational curves in ${\cal PT}$ with normal bundle ${\cal
  O}^{A}(1):={\cal O}(1)\oplus {\cal O}(1)$.  The normal bundle to
$l_x$ consists of vectors tangent to $x$ (horizontally lifted to
$T_{(x,\lambda)}{\cal F}$) modulo the twistor distribution. Therefore
we have a sequence of sheaves over $\CP^1$
\[
0\longrightarrow D \longrightarrow \C^4 \longrightarrow
{\cal O}^A(1)\longrightarrow 0.
\]
The map $\C^4 \longrightarrow {\cal O}^A(1)$ is given by
$V^{AA'}\longrightarrow V^{AA'}\pi_{A'}$.  Its kernel consists of
vectors of the form $\pi^{A'}\lambda^A$ with $\lambda^A$ varying. The
twistor distribution is therefore $D=O(-1)\otimes S^{A}$ and so there
is a canonical $L_A\in\Gamma(D\otimes {\cal O}(1)\otimes S_{A})$,
as given in (\ref{laxpair}).

%%%%%%%%%%%%%%%%%%%%%%%%%%%%%%%%%%%%%%%%%%%%%%%%%%%%%%%%%%%%%%%%%%%%%%%%%%%%%%%

\subsection{Some formulations of the ASD vacuum  condition} 
The ASD vacuum conditions  $C_{A'B'C'D'}=0$, $\Phi_{ABA'B'}=0=R$
imply the existence of a normalised, covariantly constant frame
$(o^{A'},\iota^{A'})$ of $S^{A'}$, so that $\G_{AA'B'C'}=0$.
 One can further choose an unprimed spin frame so that the Lax pair 
(\ref{laxpair}) consists of 
volume-preserving vector fields on ${\cal M}$:
\begin{prop}[Mason \& Newman \cite{MN89}.]
\label{mn}Let 
$\widehat{\nabla}_{AA'}=(\widehat{\nabla}_{00'},  \widehat{\nabla}_{01'}
\widehat{\nabla}_{10'} \widehat{\nabla}_{11'})$ 
be four independent holomorphic vector
fields on a four-dimensional complex manifold ${\cal M}$ and let $\nu$
be a nonzero holomorphic four-form. Put 
\be
\label{sde}
L_0=\widehat{\nabla}_{00'}-\l\widehat{\nabla}_{01'},
\;\;\;\;\;\; L_1=\widehat{\nabla}_{10'}-\l\widehat{\nabla}_{11'}.
\ee
Suppose that for every $\lambda \in \CP^1$
\be
\label{podstawka}
[L_0, L_1]=0,\;\;\;\;\;\;\;\;\;{\cal L}_{L_A}{\nu}=0.
\ee
Here ${\cal L}_V$ denotes the Lie derivative.
Then \[
\p_{AA'} = f^{-1}\widehat{\nabla}_{AA'},\qquad
\mbox{where}\qquad 
f^2:={\nu}(\widehat{\nabla}_{00'},  \widehat{\nabla}_{01'}
\widehat{\nabla}_{10'} \widehat{\nabla}_{11'}),
\] 
is a null-tetrad for an ASD vacuum metric. Every such metric locally
arises in this way.
\end{prop} 
In \cite{D98} the last proposition is generalised to the
hyper-Hermitian case.  A choice of unprimed spin frame with $f^2=1$ is
always possible and we shall assume this here-on so that
$\nabla_{AA'}=\widehat{\nabla}_{AA'}$.  For easy reference we rewrite the
field equations (\ref{podstawka}) in full \be
\label{eq1}
[\nabla_{A0'}, \nabla_{B0'}]=0,
\ee
\be
\label{eq2}
[\nabla_{A0'}, \nabla_{B1'}]+[\nabla_{A1'}, \nabla_{B0'}]=0,
\ee
\be
\label{eq3}
[\nabla_{A1'}, \nabla_{B1'}]=0.
\ee

Let $\Sm^{A'B'}$ be the usual basis of SD two-forms.
On the correspondence space, define 
\be
\label{PGform}
\Sm(\l):=\Sm^{A'B'}\pi_{A'}\pi_{B'}.
\ee
The formulation of the ASDVE condition dual to (\ref{podstawka}) is:
\begin{prop}[Pleba{\'n}ski \cite{Pl75}, Gindikin \cite{Gin}]
If a two-form of the form
\[\Sm(\l):=\Sm^{A'B'}\pi_{A'}\pi_{B'}\] on the correspondence space satisfies
\be
\label{Gyndykin}
\d_h \Sm(\l)=0,\;\;\;\;\;\;\;\;\;\;
\Sm(\l)\wedge\Sm(\l)=0
\ee
where $\d_h$ is the exterior derivative holding $\pi_{A'}$ constant,
then there exist one-forms $e^{AA'}$ related to $\Sm^{A'B'}$ by
equation (\ref{sdforms}) which give an ASD vacuum tetrad.
\end{prop}
Note that the simplicity condition in (\ref{Gyndykin})
arises from the condition that $\Sm^{A'B'}$ comes from a tetrad.
%Here $\d_h$ is a horizontal lift of $\d$ to $\cal F$ and so
%$\lambda$ is regarded as a parameter and is not differentiated.

To construct Gindikin's two-form starting from the twistor space, one
can pull back the fibrewise complex symplectic structure on ${\cal
  PT}\longrightarrow\CP^1$ to the projective spin bundle and fix the
ambiguity by requiring that it annihilates vectors tangent to the fibres. The
resulting two-form is $\O(2)$ valued. (To obtain Gindikin's two-form
one should divide it by a constant section of $\O(2)$.)

Put
$
\Sm^{0'0'}=-\at,\; 
\Sm^{0'1'}=\om, \;
\Sm^{1'1'}=\a\;
$. 
The second equation in (\ref{Gyndykin}) becomes
\[
\om \wedge \om =2\a\wedge\at:=-2\nu,\;\;\;
\a\wedge\om=\at\wedge\om=\a\wedge\a=\at\wedge\at =0.
\]
Equations (\ref{Gyndykin}) can be seen to arise from (\ref{podstawka}) by
observing that $\Sm(\lambda)$ can be defined by
\[
\varepsilon_{AB}\Sm(\l)=\nu(L_A, L_B, ..., ...).
\]
Note also that $L_A$ spans a two-dimensional distribution annihilating
$\Sm(\l)$. 

The two one-forms $e^A:=\pi_{A'}e^{AA'}$ by definition annihilate the
twistor distribution. Define $(1, 1)$ tensors
$\p^{B'}_{A'}:=e^{AB'}\otimes\nabla_{AA'}$ so that
\[
e^A\otimes 
L_A=\pi_{B'}\pi^{A'}\p^{B'}_{A'}=\Do +\lambda (\p -\pt )-{\lambda}^2\dd
\]
where $(\p^{0'}_{0'}, \p^{0'}_{1'}, \p^{1'}_{0'}, \p^{1'}_{1'})=
(\pt, \partial_0, \dd, \p)$. 
If the field equations are satisfied
then  the Euclidean slice of $\cal M$  is equipped with three integrable
complex structures given by $J_i:=\{i(\dd - \partial_0),\;(\p-\pt),\; 
(\dd + \partial_0)\}$
and three symplectic structures $\om_i=\{(i(\a-\at), i\om, (\a+\at)\}$
compatible with the $J_i$. It is therefore a hyper-K\"ahler manifold. 

%%%%%%%%%%%%%%%%%%%%%%%%%%%%%%%%%%%%%%%%%%%%%%%%%%%%%%%%%%%%%%%%%%%%%%%%
\subsection{The ASD condition and heavenly equations} 
Part of the residual gauge freedom in (\ref{podstawka}) is fixed by
selecting one of Pleba\'nski's null coordinate systems.
\def\wt{\tilde w}
\def\zt{\tilde z}
\begin{enumerate} 
\item
Equations (\ref{eq2}) and (\ref{eq3}) imply the existence of a coordinate system 
\[(w,z,\wt,\zt)=:(w^A,\wt^A)\] and a complex-valued function $\Om$ such that
\be
\label{oframe}
\p_{AA'}=\OM=
\Big(\frac{\p^2 \Omega}{\p w^A\p \wt^B}\frac{\p}{\p\wt_B} 
\;\;\;\frac{\p}{\p w^A} \Big ). 
\ee
Equation (\ref{eq1}) yields the first heavenly equation
\be
\label{firsteq}
\Om _{w\zt}\Om _{z\wt} -\Om _{w\wt}\Om _{z\zt}=1\;\;\mbox{or}\;\;
\frac{1}{2}\frac{\p^2\Om}{\p w_A\p \wt_B}\frac{\p^2\Om}{\p w^A\p \wt^B}=1.
\ee 
The dual tetrad is
\be
e^{A1'}=\d w^A,\;\;e^{A0'}=\frac{\p^2\Om}{\p w_A\p \wt_B}\d\wt_B
\ee
with the  flat solution  $\Om=w^A\wt_A$.
The only nontrivial part of $\Sm^{A'B'}$ is $\Sm^{0'1'}=\p\pt\Om$ so
that $\Om$ is a K\"ahler scalar. The Lax pair for the first heavenly
equation is
\begin{eqnarray}
\label{LaxH}
L_0:&=&\Om_{w\tw}\p_{\tz}-\Om_{w\tz}\p_{\tw}-\l\p_w,\nonumber\\
L_1:&=&\Om_{z\tw}\p_{\tz}-\Om_{z\tz}\p_{\tw}-\l\p_z.
\end{eqnarray}
Equations $L_0\Psi=L_1\Psi=0$ have solutions provided that
$\Om$ satisfies the
first heavenly equation (\ref{firsteq}). Here $\Psi$ is a function on 
${\cal F}$.

\item Alternatively equations (\ref{eq1}) and (\ref{eq2}) imply the
  existence of a complex-valued function $\Th$ and coordinate system
  $(w,z,x,y)=:(w^A,x_A)$, $w^A$ as above, such that \be
\label{tframe}
\p_{AA'}=\TH= 
\left( \frac{\p}{\p x^A}\;\;\; \frac{\p}{\p w^A}+
\frac{\p^2\Th}{\p x^A\p x^B}\frac{\p}{\p x_B}
\right ).
\ee
As a consequence of (\ref{eq3}) $\Th$ satisfies second heavenly equation
\be
\label{secondeq}
\Th _{xw} +\Th _{yz} + \Th _{xx}\Th _{yy}-{\Th _{xy}}^2=0\;\;\mbox{or}\;\;
\frac{\p^2\Th}{\p w^A\p x_A}+
\frac{1}{2}\frac{\p^2\Th}{\p x^B\p x^A}\frac{\p^2\Th}{\p x_B\p x_A}=0.
\ee
The dual frame is given by
\be
\label{tetr2}
e^{A0'}=\d x^A+\frac{\p^2 \Th}{\p x^B \p x_A}\d w^B,\;\; e^{A1'}=\d w^A
\ee
with $\Th=0$ defining the flat metric. The Lax pair corresponding to 
(\ref{secondeq}) is
\begin{eqnarray}
\label{Lax2}
L_0&=&\p_y-\l(\p_w-\Th_{xy}\p_y+\Th_{yy} \p_x),\nonumber\\
L_1&=&\p_x+\l(\p_z+\Th_{xx}\p_y-\Th_{xy} \p_x).
\end{eqnarray}
\end{enumerate}
Both heavenly equations were originally derived by Pleba\'nski
\cite{Pl75} from the formulation (\ref{Gyndykin}). The closure
condition is used, via Darboux's theorem, to introduce $\om^A$,
canonical coordinates on the spin bundle, holomorphic around
$\lambda=0$ such that the two-form (\ref{PGform}) is
$\Sm(\lambda)=\d_h\om^{A}\wedge \d_h\om_A$.  The various forms of the
heavenly equations can be obtained by adapting different coordinates
and gauges to these forms.

%%%%%%%%%%%%%%%%%%%%%%%%%%%%%%%%%%%%%%%%%%%%%%%%%%%%%%%%%%%%%%%%%%%%%%%%%%%%%%
\section{The recursion operator}
\label{Heavenly_Hierarchies}
\def\wt{\tilde w} \def\zt{\tilde z} 
In \S\S\ref{sectrecur} the recursion operator $R$ for the
anti-self-dual Einstein vacuum equations is constructed.  In
\S\S\ref{PTrecur} then show that the generating function for $R^i\phi$
is automatically a twistor function, and is in fact a $\check{C}$ech
representative for $\phi$.  It is shown that $R$ acts on such a
twistor function by multiplication.  A similar application to the
coordinates used in the heavenly equations yields the coordinate
description of the twistor space starting.  In \S\S\ref{methodd} we
show how that the action of the recursion operator on space-time
corresponds to multiplication of the corresponding twistor functions
by $\l$.  In \S\S\ref{algebraa} the algebra of hidden symmetries of
the second heavenly equation is constructed by applying the recursion
operator to the explicit symmetries.  In \S\S\ref{killspinors}, $R$ is
used to build a higher valence Killing spinors corresponding to hidden
symmetries. In the last subsections examples of the use of the
recursion operator are given.

%%%%%%%%%%%%%%%%%%%%%%%%%%%%%%%%%%%%%%%%%%%%%%%%%%%%%%%%%%%%%%%%%%%%%%%%%%%%%
\subsection{The recursion relations}
\label{sectrecur}
The recursion operator $R$ is a map from the space of linearised
perturbations of the ASDVE equations to itself.  This can be used to
construct the ASDVE hierarchy whose higher flows are generated by
acting on one of the coordinate flows with the recursion operator $R$.

We will identify the space of linearised perturbations
to the ASDVE equations with solutions to the background coupled wave
equations in two ways as follows.
\begin{lemma}
\label{Curvedwave}
Let $\square_{\Om}$ and $\square_{\Th}$ denote wave operators on the ASD
background determined by $\Om$ and $\Th$ respectively.
Linearised solutions to {\em(\ref{firsteq})} and {\em(\ref{secondeq})} satisfy
\be
\square_{\Om}\dom=0,\;\;\;\;\;\; \square_{\Th}\ddot=0.
\ee
\end{lemma}
{\bf Proof.} 
In both cases $\square_g=\nabla_{A1'}{\nabla^A}_{0'}$ since
\[
\square_g =\frac{1}{\sqrt{g}}\p _{a}(g^{ab}\sqrt{g}\p _{b} )=
g^{ab}\p _{a}\p _{b} + (\p _{a}g^{ab})\p _{b}
\]
but $\p _{a}g^{ab}=0$ for both heavenly coordinate systems.
For the first equation $(\p \pt (\Om +\dom ))^2 =\nu$
implies 
\[
0=(\p \pt \Om\wedge \p \pt ) \dom = \d( \p \pt \Om\wedge (\p -\pt)\dom )=
\d\ast \d\dom.
\]
Here $\ast$ is the Hodge star operator corresponding to $g$.
For the second equation we  make use of the tetrad (\ref{tframe}) and
perform coordinate calculations.\koniec

From now on we identify tangent spaces to the spaces of solutions to
(\ref{firsteq}) and (\ref{secondeq}) with the space of solutions to
the curved background wave equation, ${\cal W}_g$.  We will define the
recursion operator on the space ${\cal W}_g$.

The above lemma shows that we can consider a linearised perturbation
as an element of ${\cal W}_g$ in two ways.  These two will be related
by the square of the recursion operator.
The linearised vacuum metrics corresponding to $\dom$ and $\ddot$ are
\[
 {h ^{I}}_{AA'BB'} = \iota _{(A'}o_{B')}\nabla _{(A1'}\nabla_{B)0'}\dom,
\;\;\;\;\;
 {h ^{II}}_{AA'BB'} = o_{A'}o_{B'}\nabla _{A0'} \nabla _{B0'}\ddot.
\]
where $o^{A'}=(1,0)$ 
and  $\iota^{A'}=(0, 1)$ 
are the  constant spin frame associated to the null tetrads given above.
%We are now able to generate new linearised solutions from the old ones.
Given $\phi\in {\cal W}_g$ we use the first of these equations to find
$h^I$. If we put the perturbation obtained in this way on the LHS of the second
equation and add an appropriate gauge term  
we obtain $\phi'$ - the new element of ${\cal W}_g$ that provides
the $\ddot$ which gives rise to
\be 
\label{lmetr}
{h_{ab} ^{II}}=h_{ab} ^{I}+\nabla_{(a}V_{b)}.
\ee
To extract the recursion relations we must find $V$ such that
%$(\p -\pt )\odot \p \dom ={\cal L}_V\p \odot \pt \dom $, or equivalently
$
{h^I}_{AA'BB'}-\nabla_{(AA'}V_{BB')}=o_{A'}o_{B'}\chi _{AB}.
$
Take $V_{BB'}=o_{B'}\nabla _{B1'}\dom$, which gives
\[
\nabla_{(AA'}V_{BB')}=-\iota_{(A'}o_{B')}\nabla _{(A0'}\nabla_{B)1'}\dom 
+ o_{A'}o_{B'}\nabla _{A1'}\nabla_{B1'}\dom.
\]
This reduces ($\ref{lmetr}$) to
\be
\label{zaleznosc}
\nabla_{A1'}\nabla_{B1'}\phi =\nabla_{A0'}\nabla_{B0'}\phi'.
\ee

\begin{defi}
Define the recursion operator $R:{\cal W}_g\longrightarrow {\cal W}_g$ by
\be
\label{def}
\iota^{A'}\nabla_{AA'}\phi=o^{A'}\nabla_{AA'}R\phi,
\ee
so  formally  $R=({\nabla_{A0'}})^{-1}\circ\nabla_{A1'}$ {\em(}no summation
over the index $A${\em)}.
\end{defi}

\noindent
{\bf Remarks:}
\begin{itemize}
\item From (\ref{def}) and from 
(\ref{podstawka}) it follows that if $\phi$ belongs to ${\cal W}_g$ then
so does $R\phi$. 
\item If $R^2\dom =\ddot$ then $\dom$ and $\ddot$ correspond
to the same variation in the metric up to gauge.  
\item The operator $\phi \mapsto
\nabla_{A0'}\phi$ is 
over-determined, and its consistency follows from the wave equation on
$\phi$.  
\item This definition is formal in that in order to
invert the operator $\phi \mapsto
\nabla_{A0'}\phi$ we need to specify boundary conditions. 
\end{itemize}
To summarize: 
\begin{prop}
\label{Rekurencja}
Let ${\cal W}_g$ be the space of solutions of the wave equation on the
curved ASD background given by $g$.
\begin{enumerate}
\item[(i)] Elements of  ${\cal W}_g$ can be identified with linearised
perturbations of the heavenly equations.
\item[(ii)] There exists a (formal) 
map $R:{\cal W}_g\longrightarrow {\cal W}_g$
given by {\em(\ref{def})}.
% which generates new elements of  ${\cal W}_g$ from old.
\end{enumerate}
\end{prop}

The recursion operator can be generalised to act on solutions to the
higher helicity Zero Rest-Mass equations on the ASD vacuum backgrounds
\cite{DM97c} by using Herz potentials.  We restrict ourselves to the
gauge invariant case of left-handed neutrino field $\psi_A$ on a
heavenly background.  First note that any solution of
\[
\nabla^{AA'}\psi_{A}=0
\]
must be of the form $\nabla_{A0'}\phi$ where 
$\phi \in {\cal W}_g$. Define the recursion relations
\be
\label{conjugate}
{\cal R}\psi_A:=\nabla_{A0'}R\phi\, .
\ee
It is easy to see that ${\cal R}$ maps solutions into solutions,
although again the definition is formal in that boundary conditions
are required to eliminate the ambiguities. A conjugate recursion operator ${\cal R}$ will play a role in the  Hamiltonian formulation in Section
\ref{hformalism}. 

%%%%%%%%%%%%%%%%%%%%%%%%%%%%%%%%%%%%%%%%%%%%%%%%%%%%%%%%%%%%%%%%%%%%%%%%%%%%

\subsection{The recursion operator  and twistor functions}
\label{PTrecur}
A twistor function $f$ can be pulled back to the correspondence space
$F$.  A function $f$ on $F$ descends to twistor space iff $L_Af=0$.

Given $\phi\in{\cal W}_g$, define, for $i\in {\Z}$, a hierarchy of
linear fields, $\phi_i\equiv R^i\phi_0$. Put
$\Psi=\sum_{-\infty}^{\infty}\phi_i\lambda^i$ and observe that the
recursion equations are equivalent to $L_A\Psi=0$. Thus $\Psi$ is a
function on the twistor space $\cal PT$.  Conversely every solution of
$L_A\Psi=0$ defined on a neighbourhood of $|\lambda|=1$ can be
expanded in a Laurent series in $\lambda$ with the coefficients
forming a series of elements of ${\cal W}_g$ related by the recursion
operator.  The function $\Psi$, when multiplied by
$1/(\pi_{0'}\pi_{1'})$, is a $\check{C}$ech representative of the
element of $H^1({\cal PT, O}(-2))$ that corresponds to the solution of
the wave equation $\phi$ under the Penrose transform (i.e.\ by
integration around $|\l|=1$).  The ambiguity in the inversion of
$\nabla_{A0'}$ means that there are many such functions $\Psi$ that
can be obtained from a given $\phi$.  However, they are all equivalent
as cohomology classes.

It is clear that a series corresponding to $R\phi$ is the function
$\lambda^{-1}\Psi$.  As noted before, $R$ is not completely well defined when
acting on ${\cal W}_g$ because of the ambiguity in the inversion of
$\nabla_{A0'}$.  However, the definition $R\Psi=\Psi/\l$ is well
defined as a twistor function on $\cal PT$, but the problem resurfaces when 
one attempts to treat $\Psi(\lambda)$ as a representative of a cohomology class
since  pure gauge elements of the first sheaf
cohomology group $H^1({\cal PT, O}(-2))$ are mapped to functions defining a
non-trivial element of the cohomology.  Note,
however, that with the definition $R\Psi=\Psi/\l$, the action of $R$
is well defined on twistor functions and can be iterated without ambiguity.

%By iterating $R$ on functions and then 
%taking the corresponding cohomology classes
%we generate an infinite sequence of elements of
%H^1({\cal PT, O}(-2))$ belonging to different classes.

We can in this way build coordinate charts on twistor space from those
on space-time arising from the choices in the Plebanski reductions.
Put $\om^A_0=w^A=(w, z)$; the surfaces of constant $\om^A_0$ are
twistor surfaces. We have that ${\nabla^A}_{0'} \om^B_0=0$ so that in
particular $\nabla_{A1'}{\nabla^A}_{0'} \om^B_0=0$ and if we define
$\om^A_i=R^i\om^A_0$ then we can choose $\om^A_i=0$ for negative $i$.
% $\om^{B}=\sum_i \om^A_i\lambda^i$.
We define 
\be
\label{omegaAA}
\om ^A=\sum_{i=0}^{\infty}\om^A_i\lambda^{i}.  
\ee 
We can similarly define $\tilde{\om}^A$ by $\tilde{\om}^A_0=\wt^A$ and
choose $\tilde{\om}^A_i=0$ for $i >0$.  Note that $\om^A$ and
$\tilde{\om}^A$ are solutions of $L_A$ holomorphic around $\lambda=0$
and $\lambda=\infty$ respectively and they can be chosen so that they
extend to a neighbourhood of the unit disc and a neighbourhood of the
complement of the unit disc and can therefore be used to provide a patching
description of the twistor space.

%%%%%%%%%%%%%%%%%%%%%%%%%%%%%%%%%%%%%%%%%%%%%%%%%%%%%%%%%%%%%%%%%%%%%%%%%%%%%%
\subsection{The Penrose transform of 
linearised deformations and the recursion operator}
\label{methodd} 
The recursion operator acts on linearised perturbations of the ASDVE
equations.  Under the twistor correspondence, these correspond to
linearised holomorphic deformations of (part of) ${\cal PT}$.

Cover $\cal PT$ by two sets, $U$ and $\tilde{U}$ with
$|\lambda|< 1+\epsilon$ on $U$ and $|\lambda|>1-\epsilon$ on
$\tilde{U}$ with $(\om^A,\lambda)$ coordinates on $U$ and
$(\tilde{\om}^A,\lambda^{-1})$ on $\tilde{U}$. The twistor space 
$\cal PT$ is then
determined by the transition function ${\tilde{\om}}^B=
{\tilde{\om}}^B(\om^A, \pi_{A'})$ on $U\cap \tilde{U}$ which preserves
the fibrewise 2-form, $\d
\omega^A\wedge\d\omega_A|_{\l=\mathrm{const.}} =
  \d\tilde{\omega}^A\wedge\d\tilde{\omega}_A|_{\l=\mathrm{const.}}$. 

Infinitesimal deformations are given by elements of   
$H^1({\cal PT}, {\bf \Th})$, where ${\bf \Th}$ denotes a sheaf of
germs of holomorphic vector fields. Let 
\[
Y=f^A(\om^B, \pi_{B'})\frac{\p}{\p \om^A}
\]
defined on the overlap $U\cap \tilde{U}$ and define a class in $
H^1({\cal PT}, {\bf \Th})$ that preserves the fibration ${\cal
  PT}\mapsto \CP^1$. The corresponding infinitesimal deformation
is given by
\be
\label{ideform}
\tilde{\om}^A(\omega^A,\pi_{A'},t)=(1+tY)(\tilde{\om}^A)+O(t^2).
\ee
From the globality of $\Sm(\lambda)=\d\om^A\wedge \d\om_A$ it follows that $Y$
is a Hamiltonian vector field with a Hamiltonian
$f\in H^1({\cal PT,O}(2))$ with respect to the symplectic structure $\Sm$.
A finite deformation is given by integrating
\[
\frac{\d\tilde{\om}^B}{\d t}=\varepsilon^{BA}\frac{\p f}{\p\tilde{\om}^A}. 
\]
from $t=0$ to $1$. 
Infinitesimally we can put
\be
\delta\tilde{\om}^A=\frac{\p\delta f}{\p\tilde{\om}_A} .
\ee
%This can be understood as follows:
%$\tilde{\om}^A$ is the patching
%function obtained by exponentiating the Hamiltonian vector field of
%$f$ and corresponds to
If the ASD metric is determined by $\Th$ and then
$\varepsilon^{BA}\p\delta f/\p \omega^B$,  (or more simply $\delta f$)
is a linearised deformation corresponding to $\delta \Th\in {\cal W}_g$. 

The recursion operator acts on linearised deformations as follows
\begin{prop}
\label{Twierdzenie1}
Let $R$ be the recursion operator defined by {\em(\ref{def})}. Its twistor 
counterpart is the multiplication operator
\be
\label{twierdzenie}
R\;\delta f=\frac{\pi_{1'}}{ \pi_{0'}}\delta f=\lambda^{-1}\delta f.
\ee
\end{prop}
%We see that the twistor description of the recursion operator
%is simpler than the space-time one. 

\noindent
[Note that $R$ acts on $\delta f$ without ambiguity; the ambiguity in
boundary condition for the definition of $R$ on space-time is absorbed
into the choice of explicit representative for the cohomology class
determined by $\delta f$.]

\smallskip

\noindent 
{\bf Proof.} 
Pull back  $\delta f$ to the primed spin bundle on which it is a
coboundary so that
\be
\label{split}
\delta f(\pi_{A'}, x^a)=h(\pi_{A'}, x^a)-\tilde{h}(\pi_{A'}, x^a) \ee
where $h$ and $\tilde{h}$ are holomorphic on $U$ and $\tilde{U}$
respectively (here we abuse notation and denote by $U$ and $\tilde{U}$
the open sets on the spin bundle that are the preimage of $U$ and
$\tilde{U}$ on twistor space). A choice for the splitting
(\ref{split}) is given by
\begin{eqnarray}
\label{hl}
h&=&\frac{1}{2\pi i}\oint_{\Gamma}
\frac{(\pi^{A'}o_{A'})^3}{(\rho^{C'}\pi_{C'}) (\rho^{B'}o_{B'})^3}
\delta f(\rho_{E'})\rho_{D'}\d\rho^{D'},\\
\tilde{h}&=&\frac{1}{2\pi i}\oint_{\tilde{\Gamma}}
\frac{(\pi^{A'}o_{A'})^3}{(\rho^{C'}\pi_{C'}) (\rho^{B'}o_{B'})^3}
\delta f(\rho_{E'})\rho_{D'}\d\rho^{D'}.\nonumber
\end{eqnarray}
Here $\rho_{A'}$ are homogeneous coordinates of $\CP^1$ pulled back to
the spin bundle.  The contours ${\Gamma}$ and $\tilde{\Gamma}$ are
homologous to the equator of $\CP^1$ in $U\cap\tilde{U}$ and are such
that ${\Gamma}-{\tilde{\Gamma}}$ surrounds the point
$\rho_{A'}=\pi_{A'}$.

The functions $h$ and $\tilde{h}$ are homogeneous of degree 1 in
$\pi_{A'}$ and do not descend to ${\cal PT}$, whereas their difference
does so that \be
\label{potential}
\pi^{A'}\nabla_{AA'}h=\pi^{A'}\nabla_{AA'}\tilde{h}=
\pi^{A'}\pi^{B'}\pi^{C'}\Sm_{AA'B'C'} \ee where the first equality
shows that the LHS is global with homogeneity degree 2 and implies the
second equality for some $\Sm_{AA'B'C'}$ which will be the third
potential for a linearised ASD Weyl spinor.  $\Sm_{AA'B'C'}$ is in
general defined modulo terms of the form $\nabla_{A(A'}\gamma_{B'C')}$
but this gauge freedom is partially fixed by choosing the integral
representation above; $h$ vanishes to third order at $\pi_{A'}=o_{A'}$
and direct differentiation, using $\nabla_{AA'}\delta f=\rho_{A'}\delta
f_A$ for some $\delta f_A$, gives
$\Sm_{AA'B'C'}=o_{A'}o_{B'}o_{C'}\nabla_{A0'}\ddot$ where \be
\label{contur}
\ddot=\frac{1}{2\pi i} \oint_{\Gamma}\frac{\delta
  f}{{(\rho^{B'}o_{B'})}^4} \rho_{D'}\d\rho^{D'}.  
\ee 
This is consistent with the Plebanski gauge choices (there is also a
gauge freedom in $\ddot$ arising from cohomology freedom in $\delta f$
which we shall describe in the next subsection.)  The condition
$\nabla_{A(D'}{\Sm^A}_{A'B'C')}=0$ follows from equation
(\ref{potential}) which, with the Pleba\'nski gauge choice, implies
$\ddot\in {\cal W}_g$.  Thus we obtain a twistor integral formula for
the linearisation of the second heavenly equation.

%(\ref{potential}) becomes
%\[
%\pi^{A'}\nabla_{AA'}\frac{1}{2\pi i}\oint_{\Gamma}
%\frac{\delta f(\rho_{E'})}{(\rho^{C'}\pi_{C'}) (\rho^{B'}o_{B'})^3}
%\rho_{D'}d\rho^{D'}=\nabla_{A0'}\dot.
%\]
%Define $\delta f_A$ by $\nabla_{AA'}\delta f=\rho_{A'}\delta f_A$.
%Equation (\ref{potential}) becomes
%\be
%\label{fa}
%\oint_{\Gamma}\frac{\delta f_A(\rho_{E'})}{(\rho^{B'}o_{B'})^3}
%\rho_{D'}\d\rho^{D'}=2\pi i\nabla_{A0'}\ddot.
%\ee
%The twistor function  ${\delta f}$ is not constrained by the 
%RHS of (\ref{fa})
%being a gradient. To see this define $\delta f_{AB}$ by  
%$\nabla_{AA'}(\delta f_B \rho_{B'})=\delta f_{AB}\rho_{A'}\rho_{B'}$
%and note that in the ASD vacuum $\delta f_{AB}$ is symmetric which implies
%${\nabla^A}_{A'}\delta f_A=0$. Therefore the RHS of   
%(\ref{fa}) is also a solution of a neutrino equation so (in the ASD vacuum)
%it must be given by $\a^{A'}\nabla_{AA'}\phi$ where $\a^{A'}$ is a constant
%spinor and $\phi \in {\cal W}_g$. Equation (\ref{fa}) gives the formula 
%for a linearisation of the second heavenly equation

Now recall formula (\ref{def}) defining $R$. Let $R\delta f$ be the
twistor function corresponding to $R\ddot$ by (\ref{contur}). The
recursion relations yield
\[
\oint_{\Gamma}\frac{R\delta f_A}{(\rho^{B'}o_{B'})^3}
\rho_{D'}\d\rho^{D'}=\oint_{\Gamma}\frac
{\delta f_A}{(\rho^{B'}o_{B'})^2(\rho^{B'}\iota_{B'})}
\rho_{D'}\d\rho^{D'}
\]
so $R\delta f=\lambda^{-1}\delta f$. \koniec

\smallskip

Let $\dom$ be the linearisation of the first heavenly potential.
From $R^2\dom=\ddot$ it follows that
\[
\dom=\frac{1}{2\pi i}\oint_{\Gamma}\frac{\delta f}{(\rho_{A'}o^{A'})^2
(\rho_{B'}\iota^{B'})^2}
\rho_{C'}\d\rho^{C'}.
\]
%%%%%%%%%%%%%%%%%%%%%%%%%%%%%%%%%%%%%%%%%%%%%%%%%%%%%%%%%%%%%%%%%%%%%%%%%%%%%%
\subsection{Hidden symmetry algebra}
\label{algebraa}
The ASDVE equations in the Pleba\'nski forms have a residual coordinate
symmetry.  This consists of area preserving diffeomorphisms in the
$w^A$ coordinates together with some extra transformations that depend
on whether one is reducing to the first or second form.  By
regarding the infinitesimal forms of these transformations as
linearised perturbations and acting on them using the recursion
operator, the coordinate (passive) symmetries can be extended to give 
`hidden' (active) symmetries of the heavenly equations.  Formulae
(\ref{contur}) and (\ref{twierdzenie}) can be used to recover the
known relations (see for example \cite{Ta89}) of the hidden symmetry
algebra of the heavenly equations.  We deal with the second equation
as the case of the first equation was investigated by other methods
\cite{Ph90}.

Let $M$ be a volume preserving vector field on $\cal M$. Define
$\delta_M^0\nabla_{AA'}:=[M, \nabla_{AA'}]$.  This is a pure gauge
transformation corresponding to addition of ${\cal L}_Mg$ to the
space-time metric and preserves the field equations. Note that
\[
[\delta_M^0, \delta_N^0]\nabla_{AA'}:=\delta_{[M, N]}^0\nabla_{AA'}.
\]
Once a Pleba\'nski coordinate system and reduced equations have been
obtained, the reduced equation will not be invariant under all the
SDiff$(\cal M)$ transformations.  The second form will be preserved if
we restrict ourselves to transformations which preserve the SD
two-forms $\Sm^{1'1'}=\d w_A\wedge \d w^A$ and $\Sm^{0'1'}=\d
x_A\wedge \d w^A$. The conditions ${\cal L}_M\Sm^{0'0'}={\cal
  L}_M\Sm^{0'1'}=0$ imply that $M$ is given by
\[
 M= \frac{\p h}{\p w_A }  \frac{\p}{\p w^A} +\Big(\frac{\p g}{\p w_A}-
x^B\frac{\p^2 h}{\p w_A \p w^B} \Big)\frac{\p}{\p x^A}
\]
where $h=h(w^A)$ and $g=g(w^A)$. The space-time is now viewed as a
cotangent bundle ${\cal M}=T^*{\cal N}^2$ with $w^A$ being coordinates
on a two-dimensional complex manifold ${\cal N}^2$. The full
SDiff$({\cal M})$ symmetry breaks down to the semi-direct product of
SDiff$({\cal N}^2)$, which acts on ${\cal M}$ by a Lie lift, with
$\Gamma ( {\cal N}^2, {\cal O})$ which acts on ${\cal M}$ by
translations of the zero section by the exterior derivatives of
functions on ${\cal N}^2$. Let $\delta_M\Th$ correspond to
$\delta_M^0\nabla_{AA'}$ by
\[
\delta_M^0\nabla_{A1'}=\frac{\p^2\delta_M\Th}{\p x^A\p
x^B}\frac{\p}{\p x_B}.
\]
The `pure gauge' elements  are  
\begin{eqnarray}
\label{trtr}
\delta^0_M \Th &=&F +x_AG^A+
x_Ax_B\frac{\p^2 g}{\p w_A \p w_B}+x_Ax_Bx_C \frac{\p^3 h}{\p w_A \p w_B \p
w_C}\nonumber\\
& &+ \frac{\p g}{\p w_A }\frac{\p \Th}{\p x^A }+
\frac{\p h}{\p w_A }\frac{\p \Th}{\p w^A }
-x^B\frac{\p^2 h}{\p w_A \p w^B}\frac{\p \Th}{\p x^A }
\end{eqnarray}
where $F, G^A, g, h$ are functions of $w^B$ only.

The above symmetries can be seen to arise from symmetries on twistor
space as follows.  Since we have the symplectic form $\Sm=\d\om^A\wedge
\d\om_A$ on the fibres of $\mu:{\cal PT}\longrightarrow \CP^1$, a
symmetry is a holomorphic diffeomorphism of the set $U$ that restricts to a 
canonical transformation on each fibre.
%of $\mu$ leaving $\Sm$ invariant on a neighbourhood of $\lambda=0$.  
Let $H=H(x^a,\lambda)=\sum_{i=0}^{\infty} h_i\lambda^i$ be the
Hamiltonian for an infinitesimal such transformation pulled back to
the projective spin bundle.  The functions $h_i$ depend on space time
coordinates only. In particular $h_0$ and $h_1$ give 
$h$ and $g$ from the previous construction (\ref{trtr}).  This can be
seen by calculating how $\Th$ transforms if $\om^A=w^A+\lambda x^A
+\lambda^2 \p\Th/\p x_A+...  \longrightarrow \hat{\om}^A$. Now $\Th$
is treated as an object on the first jet bundle of a fixed fibre of
${\cal PT}$ and it determines the structure of the second jet.

These symmetries take a solution to an equivalent solution.  The
recursion operator can be used to define an algebra of `hidden
symmetries' that take one solution to a different one as follows. 

Let $\delta^0_M\Th$ be an expression of the form 
(\ref{trtr}) which also satisfies $\Box_g\delta^0_M\Th=0$.
We
set $${\delta_M}^i\Th:=R^i{\delta_M}\Th\in {\cal W}_g.$$

\begin{prop}
Generators of the hidden symmetry algebra of the second heavenly equation
satisfy the relation
\be
\label{walg}
[{\delta_M}^i, {\delta_N}^j]={\delta_{[M,N]}}^{i+j}.
\ee
\end{prop}

\noindent
{\bf Proof.}
This can be proved directly by showing that the ambiguities in $R$ can
be chosen so that $R\circ \delta_M=\delta_M\circ R$.  It is perhaps
more informative to prove it by its action on twistor functions.

Let $\delta_M^i f$ be the twistor function corresponding to
$\delta_M^i\Theta$ (by (\ref{contur})) treated as an element of
$\Gamma(U \cap \tilde{U}, {\cal O}(2))$ rather than $H^1({\cal
  PT,O}(2))$.  Define $[\delta_M^i, \delta_N^j]$ by
\[
[\delta_M^i, \delta_N^j]\Th:=\frac{1}{2\pi i}\oint
\frac{\{\delta_M^i f,\delta_N^j f\}}{(\pi_{0'})^4}\pi_{A'}\d\pi^{A'}
\]
where the Poisson bracket is calculated with respect to a canonical 
Poisson structure on ${\cal PT}$. From Proposition (\ref{twierdzenie})
it follows that
\[
[\delta_M^i, \delta_N^j]\Th=\frac{1}{2\pi i}\oint\lambda^{-i-j}
\frac{\{\delta_M f,\delta_N f\}}{(\pi_{0'})^4}\pi_{A'}\d\pi^{A'}
%=\frac{1}{2\pi i}\oint
%\frac{\Sm(Y^i_M,Y^j_N)}{(\pi_{0'})^4}\pi_{A'}d\pi^{A'}
=R^{i+j}\delta_{[M,N]}\Th
\]
as required.\koniec
\subsection{Recursion procedure for Killing spinors}
\label{killspinors}
Let $({\cal M}, g)$ be an ASD vacuum space. We say that $L_{A_1'...A_n'}$
is a Killing spinor of type $(0, n)$ if 
\be
\label{Kspinor}
{\nabla^A}_{(A'}L_{B_1'...B_n')}=0.
\ee
Killing spinors of type $(0, n)$ give rise to 
Killing spinors of type $(1, n-1)$ by
\[
{\nabla^A}_{A'}L_{B_1'...B_n'}=\varepsilon_{A'(B_1'}{K^A}_{B_2'...B_n')}.
\]
In an ASD vacuum, $K^{BB_2'...B_n'}$ is also a Killing spinor 
\[
{\nabla^{(A}}_{(A'}{K^{B)}}_{B_1'...B_n')}=0.
\]
Put (for $i=0, ..., n$)
\[
L_i:=\iota^{B_1'}...\iota^{B_i'}o^{B_{i+1}'}...o^{B_n'}L_{B_1'...B_n'},
\]
and contract (\ref{Kspinor}) with
$\iota^{B_1'}...\iota^{B_i'}o^{B_{i+1}'}...o^{B_{n+1}'}$ to obtain
\[
i\nabla_{A1'}L_{i-1}=-(n-i+1)\nabla_{A0'}L_i,\qquad i=0, ..., n-1.
\]
We  make use of  the recursion relations (\ref{def}):
\[
\frac{-i}{n+1-i}R(L_{i-1})=L_i.
\]
This leads to a general formula for Killing spinors (with
$\nabla_{A0'}L_0=0$)
\be
L_i=(-1)^i{n\choose i}^{-1}R^i(L_0),
\qquad L_{B_1'B_2'...B_n'}=\sum_{i=0}^n
o_{(B_1'}...o_{B_i'}\iota_{B_{i+1}'}...\iota_{B_n')}L_i
\ee
and equation (\ref{Kspinor}) is then satisfied iff $R^{-1}L_0=RL_n=0$.

%%%%%%%%%%%%%%%%%%%%%%%%%%%%%%%%%%%%%%%%%%%%%%%%%%%%%%%%%%%%%%%%%%%%%%%%%%%%%%%
\subsection{Example 1}  Let us demonstrate 
how to use the recursion procedure to
find metrics with hidden symmetries.
Let  $\p_{t_n}\Om:=\phi_n$ be a linearisation of the first heavenly
equation. We have $R:z\longrightarrow\Om_w=\p_{t_1}\Om$. 
Look for solutions to (\ref{firsteq})
with an additional constraint $\p_{t_2}\Om=0$. The recursion
relations (\ref{def}) imply $\Om_{wz}=\Om_{ww}=0$, therefore
\[
\Om(w, z, \tw, \tz)=wq(\tw, \tz)+P(z, \tw, \tz).
\]

The heavenly equation yields $\d q\wedge\d P\wedge \d z=\d \tz \wedge
\d \tw \wedge \d z$. With the definition $\p_z P=p$ the metric is
\[
\d s^2=2\d w\d q+2\d z\d p+f\d z^2,
\]
where $f=-2P_{zz}$. We adopt $(w, z, q, p)$ as a new coordinate system. 
Heavenly equations imply that $f=f(q, z)$ is an arbitrary function of
two variables. These are the null ASD plane wave solutions.

\subsection{Example 2}
Now  we shall illustrate the Propositions
\ref{Rekurencja} and \ref{Twierdzenie1} with the example of 
the Sparling--Tod solution \cite{ST79}.
The coordinate formulae for the pull back of twistor
functions are:
\begin{eqnarray}
\label{tfunctions}
{\mu}^0&=&w+\lambda
y-\lambda^2\Th_x+ \lambda^3\Th_z +...\;\;,\nonumber\\
{\mu}^1&=&z-\lambda x-\lambda^2\Th_y-\lambda^3\Th_w +...\;.
\end{eqnarray}
Consider  
\be
\label{stsol}
\Th=\frac{\sigma}{wx+zy},
\ee
where $\sigma=const$. It satisfies both the linear and the nonlinear
part of (\ref{secondeq}).

\noindent
{\bf The flat case:} First we shall treat (\ref{stsol}), with
$\sigma=1$, as a solution $\phi_0$ to the wave equation
on the flat background.  The recursion relations are
\[
(R\phi_0)_x=\frac{y}{(wx+zy)^2},\;\;\;(R\phi_0)_y=\frac{-x}{(wx+zy)^2}.
\]
They have a solution $\phi_1:=R\phi_0=(-y/w)\phi_0$. More generally we
find that 
\be
\label{flatex}
\phi_n:=R^n\phi_0=\Big(-\frac{y}{w}\Big)^n \frac{1}{wx+zy}.
\ee
The last formula can be also found  using twistor methods.
The twistor function corresponding to $\phi_0$ is $1/(\mu^0\mu^1)$,
where $\mu_0=w+\lambda y$ and $\mu_1=z-\lambda x$. 
By Proposition \ref{twierdzenie} the
twistor function corresponding to $\phi_n$ is $\l^{-n}/(\mu^0\mu^1)$.
This can be seen by applying the formula (\ref{contur}) and
computing the residue at the pole $\l=-w/y$. It is interesting to ask
whether any $\phi_n$ (apart from $\phi_0$) is a solution to the 
heavenly equation. Inserting $\Th=\phi_n$ to (\ref{secondeq}) yields
$n=0$ or $n=2$. We parenthetically mention that $\phi_2$ yields (by
formula (\ref{tetr2})) a metric of type $D$ which is conformal to the
Eguchi-Hanson solution.
%%%%%%%%%%%%%%%%%%%%%%%%%%%%%%%%%%%%%%%%%%%%%%%%%%%%%%%%%%%%%%%%%%%%%%%%%%%%%%%

\noindent
{\bf The curved case.} Now let $\Th$ given by (\ref{stsol}) determine the curved metric
\be
\label{SparlingTod}
\d s^2=2\d w\d x+2\d z \d y+ 4\sigma(wx+zy)^{-3}(w\d z-z\d w)^2.
\ee
The recursion relations
\[
\p_y(R\phi)=(\p_w-\Th_{xy}\p_y+\Th_{yy} \p_x)\phi,\;\;\;
-\p_x(R\phi)=(\p_z+\Th_{xx}\p_y-\Th_{xy} \p_x)\phi
\]
are
\begin{eqnarray*}
-\p_x(R\psi)&=&(\p_z+2\sigma w(wx+zy)^{-3}(w\p_x-z\p_y))\psi,
\\
\p_y(R\psi)&=&(\p_w+2\sigma z(wx+zy)^{-3}(w\p_x-z\p_y))\psi,
\end{eqnarray*}
where $\psi$ satisfies
\be
\label{stwave}
\square_{\Th}\psi=
2(\p_x\p_w+\p_y\p_z+
2\sigma(wx+zy)^{-3}(z^2{\p_x}^2+w^2{\p_y}^2-2wz\p_x \p_y))\psi
=0.
\ee
One solution to the last equation is $\psi_1=(wx+zy)^{-1}$.
We apply the recursion relations 
to find  the sequence of linearised solutions 
\begin{eqnarray*}
\psi_2&=&\Big(-\frac{y}{w}\Big)\frac{1}{wx+zy},\;
\psi_3= -\frac{2}{3}\frac{\sigma}{(wx+zy)^3}+ 
\Big(-\frac{y}{w}\Big)^2\frac{1}{wx+zy}, ...,\\
\psi_n&=&\sum_{k=0}^nA_{(n)}^k\Big(-\frac{y}{w}\Big)^k
(wx+zy)^{k-n}.
\end{eqnarray*}
To find $A_{(n)}^k$ note that the recursion relations imply
\begin{eqnarray*}
&R&\Big(\Big(-\frac{y}{w}\Big)^k(wx+zy)^{j}\Big)=\\
&=&\Big(\Big(-\frac{y}{w}\Big)
-\sigma\Big(-\frac{y}{w}\Big)^{-1}(wx+zy)^{-2}\frac{k}{j+2}\Big)
\Big(-\frac{y}{w}\Big)^k(wx+zy)^{j}\Big).
\end{eqnarray*}
This yields a recursive formula
\be
A_{(n+1)}^{k}=A_{(n)}^{k-1}-2\sigma\frac{k+1}{n-k+1}A_{(n)}^{k+1},\;\;\;\;\;
A_{(1)}^0=1,\;\;A_{(1)}^1=0,\;\;A_{(n)}^{-1}=0,\;\;k=0...n,
\ee
which determines the algebraic (as opposed to the differential) recursion
relations between $\psi_n$ and $\psi_{n+1}$.
It can be checked that functions $\psi_n$ indeed satisfy (\ref{stwave}).
Notice that if $\sigma=0$ (flat background) then we recover (\ref{flatex}). 
We can also find the inhomogeneous twistor coordinates pulled back to 
${\cal F}$
\begin{eqnarray*}
{\mu}^0&=&w+\lambda y+
\sum_{n=0}^{\infty}\sigma\l^{n+2}\sum_{k=0}^nB_{(n)}^kw\Big(-\frac{y}{w}\Big)^k
(wx+zy)^{k-n-1},\\
{\mu}^1&=&z-\lambda x+
\sum_{n=0}^{\infty}\sigma\l^{n+2}\sum_{k=0}^nB_{(n)}^kz\Big(\frac{x}{z}\Big)^k
(wx+zy)^{k-n-1}.
\end{eqnarray*}
where
\[
B_{(n+1)}^{k}=B_{(n)}^{k-1}-2\sigma\frac{k+1}{n-k+2}B_{(n)}^{k+1},\;\;\;\;\;
B_{(1)}^0=1,\;\;B_{(1)}^1=0,\;\;B_{(n)}^{-1}=0,\;\;k=0...n.
\]
The polynomials $\mu^A$ solve  $L_{A}(\mu^B)=0$,
where now
\begin{eqnarray*}
L_0&=&-\l \p_w -2\l\sigma z^2(wz+zy)^{-3}\p_x  +(1+ 2\l\sigma
wz(wz+zy)^{-3})\p_y,\\ 
L_1&=&\l\p_z +(1-2\l\sigma wz(wz+zy)^{-3})\p_x +2\l\sigma
w^2(wz+zy)^{-3})\p_y. 
\end{eqnarray*}
%%%%%%%%%%%%%%%%%%%%%%%%%%%%%%%%%%%%%%%%%%%%%%%%%%%%%%%%%%%%%%%%%%%%%%%%%%%%%%%

\section{Hierarchies for the ASD vacuum equations}

\label{Heavenly_h}

\def\wt{\tilde w} \def\zt{\tilde z} The hidden symmetries
corresponding to higher flows associated to translations along the
coordinate vector fields give `higher flows' of a
hierarchy.  This yields a hierarchy of flows of the
anti-self-dual Einstein vacuum equations.  We first give this for the
equations in their second heavenly form but then give the equations in
the form of consistency conditions for a Lax system of vector fields
generalizing equations \ref{podstawka}.  The nonlinear graviton
construction generalizes to give a construction for the corresponding
system of equations and is presented in \S\S\ref{ONONtwistor}.  In
\S\S\ref{geomstr} the geometric structure of solutions to the
truncated hierarchy are explored in further detail.  Finally in
\S\S\ref{infdef} infinitesimal deformations are studied.

%\smallskip
%; the heavenly hierarchies are shown to arise naturally on
%the moduli space of certain rational curves in a twistor space.

\subsection{Hierarchies for the heavenly equations}
The generators of higher flows are first obtained by applying powers
of the recursion operator to the linearised perturbations
corresponding to the evolution along coordinate vector fields.  This
embeds the second heavenly equation into an infinite system of
over-determined, but consistent, PDEs (which we will truncate at some
arbitrary but finite level).  These equations in turn can be naturally
embedded into a system of equations that are the consistency
conditions for an associated linear system that extends
(\ref{podstawka}).  We shall discuss here the hierarchy for the second
Pleba\'nski form; that for the first arises from a different
coordinate and gauge choice.

Introduce the coordinates $x^{Ai}$, where for $i=0, 1, x^{Ai}=x^{AA'}$
are the original coordinates on ${\cal M}$, and for $1<i \leq n,
x^{Ai}$ are the parameters for the new flows (with $2n-2$ dimensional
parameter space $\X$).  The propagation of $\Th$ along these
parameters is determined by the recursion relations
\begin{eqnarray}
\label{hierdef}
\p_y(\p_{Bi+1}\Theta) &=& (\p_w-\Th_{xy}\p_y+\Th_{yy} \p_x)
\p_{Bi}\Theta \, ,\nonumber \\
-\p_x(\p_{Bi+1}\Theta)&=&(\p_z+\Th_{xx}\p_y-\Th_{xy}
\p_x)\p_{Bi}\Theta\, ,\nonumber \\
\mbox{or } \quad \p_{A0} (\p_{Bi+1}\Theta)& =& (\p_{A1}+
\p_{C0}\p_{A0}\Theta\p^C{}_0 )\p_{Bi}\Theta\, .
\end{eqnarray}
%together with some additional conditions, given below, that partially
%fix the ambiguity in the definition of the recursion operator.
%This gives 
%the explicit field equations
%\be\label{hierdef}
%\p_{A1'}\p_{Bj-1}\Th
% +\p_{A0'}\p_{C0'}\Th{\p^C}_{0'}\p_{Bj-1}\Th= 
%\p_{A0'}\p_{Bj}\Th.
%\ee
However, we will take the hierarchy to be the system (containing the
above when $j=1$) 
%see that consistency conditions for (\ref{hierdef}) imply that in
%addition $\Th$ satisfies  the equations 
\be
\label{hier2}
\p_{Ai}\p_{Bj-1}\Th -\p_{Bj}\p_{Ai-1}\Th +\{ \p_{Ai-1}\Th, \p_{Bj-1}\Th
\}_{yx}=0,\;\;\;\; i,j=1...n.
\ee
Here $\{..., ...\}_{yx}$ is the Poisson bracket with respect to 
the Poisson structure $\p/\p x^A\wedge\p/\p x_A=2\p_x\wedge\p_y$.

\begin{lemma}
The linear system  for equations {\em(\ref{hier2})} is
\be
\label{system}
L_{Ai}s=(-\lambda D_{Ai+1}+\delta_{Ai})s=0, \qquad i=0, ..., n-1, 
\ee
where
\begin{enumerate}
\item $s:=s(x^{Ai}, \lambda )$ is a function on a spin bundle (a
$\CP^1$-bundle) over ${\cal N}={\cal M}\times \X$, 
\item
$D_{Ai+1}:=\p_{Ai+1} +[\p_{Ai}, V]$,
{\em($V=\varepsilon^{AB}\p_{A0}\Th\p_{B0}$)\em}
and $\delta_{Ai}:=\p_{Ai}$ are $4n$ vector fields on ${\cal N}$. 
\end{enumerate}
\end{lemma}
{\bf Proof.}  This follows by direct calculation.
The compatibility conditions for (\ref{system}) are: 
\be
\label{one}
[D_{Ai+1}, D_{Bj+1}]=0,
\ee
\be
\label{two}
[\delta_{Ai}, \delta_{Bj}]=0,
\ee
\be
\label{three}
[D_{Ai+1}, \delta_{Bj}]-[D_{Bj+1}, \delta_{Ai}]=0.  
\ee 
It is
straightforward to see that equations (\ref{two}) and (\ref{three})
hold identically with the above definitions and (\ref{one}) is
equivalent to (\ref{hier2}). \koniec

%Equation (\ref{one})
%then holds for $i=0$ and $j$ arbitrary as a consequence of equation
%(\ref{hierdef}).  To obtain  equations 
%(\ref{hier2})  we must assume that we
%have chosen the ambiguity in our solution to the recursion operator in
%such a way that $[D_{Ai+1},D_{Bj+1}]=0$ at some value of $x^{A0}$. 
%Then the Jacobi identity will give
%\begin{eqnarray*}
%0&=& [L_{A0}, [L_{Bi}, L_{Cj}]]= [L_{A0}, [D_{Bi+1}, D_{Cj+1}]]\\
%&=&[\p_{A0}, [D_{Bi+1}, D_{Cj+1}]]  -\lambda  [D_{A1}, [D_{Bi+1},
%D_{Cj+1}]]\, 
%\end{eqnarray*}
%so that in particular $[\p_{A0}, [D_{Bi+1}, D_{Cj+1}]] =0$.  This,
%with the assumption that $[D_{Ai+1},D_{Bj+1}]=0$ at some value of $x^{A0}$
%implies that that $[D_{Ai+1}, D_{Bj+1}]=0$ everywhere 
%(note that this will follow from weaker assumptions as we must 
%also have $[D_{A1}, [D_{Bi+1},
%D_{Cj+1}]]=0$ which implies that we only need assume 
%$[D_{Ai+1},D_{Bj+1}]=0$ at some fixed value of $(x^{A0}, x^{A1})$). The 
%equations $[D_{Ai+1},D_{Bj+1}]=0$ yield the additional equations 
%(\ref{hier2}) on $\Th$. \koniec

\smallskip

As a converse to this lemma, we will see in \S\S\ref{ONONtwistor}
using the twistor correspondence, that given the Lax system above, in
which the vector fields $D_{Ai}$ and $\delta_{Aj}$ are volume
preserving vector fields, then coordinate and gauge choices can be
made so that the Lax system takes on the above form.

\subsubsection{Spinor notation}
The above can also be represented in a spinorial formulation that
will be useful later. We introduce the spinor indexed coordinates
$x^{AA_1'...A_n'}=x^{A(A_1'...A_n')}$ on ${\cal N}$ 
which correspond to the $x^{Ai}$ by
\[
x^{Ai}={n \choose i}x^{AA_1'A_2'...A_n'}o_{A_1'}...o_{A_i'}
\iota_{A_{i+1}'}...\iota_{A_{n}'}(-1)^{n-i}.
\]
The vector fields $D_{Ai+1}$ and $\delta_{Ai}$ are then represented by
the $4n$ vector fields on ${\cal N}$, $D_{AA_1'(A_2'...A_n')}$ where
\[ 
D_{AA_1'i}=\iota^{A_2'}...\iota^{A_{i}'}o^{A_{i+1}'}...o^{A_{n}'}
D_{AA_1'(A_2'...A_n')}, \;\;D_{A1'i}=D_{Ai+1},\;\;D_{A0'i}=\delta_{Ai}
\]  
and 
$L_{A(A_2'...A_n')}=\pi^{A_1'}D_{AA_1'(A_2'...A_n')},\;\;
L_{Ai}=\pi^{A_1'}D_{AA_1'i}$.
In the adopted gauge 
\[
D_{A0'A_2'...A_n'}=\p_{A0'A_2'...A_n'},\;\;\;
D_{A1'A_2'...A_n'}=\p_{A1'A_2'...A_n'}+[\p_{A0'A_2'...A_n'},V].
\]
In what follows we will often be interested in 
$\nabla_{A(A_1'A_2'...A_n')}$, the symmetric part of 
$D_{AA_1'A_2'...A_n'}$.
\begin{eqnarray}
\nabla_{Ai}&=&D_{A(A_1'A_2'...A_n')}
\iota^{A_1'}...\iota^{A_{i}'}o^{A_{i+1}'}...o^{A_{n}'}\\
&=&\frac{1}{n}(iD_{A1'i-1}+(n-i)D_{A0'i})
=\p_{Ai}+\frac{i}{n}[\p_{Ai-1}, V].
\end{eqnarray}
Put $D_{A0'...0'}=\p_A$. The $2n+2$ vector fields \[
\nabla_{AA_1'...A_n'}=\{ \p_A, \nabla_{A0'1'A_2'...A_{n-1}'},D_{An}\}
\]
span $T^*{\cal N}$.

\subsection{The twistor space for the hierarchy}
\label{ONONtwistor}
The twistor space ${\cal PT}$ for a solution to the hierarchy
associated to the Lax system on ${\cal N}$ as above is obtained by
factoring the spin bundle ${\cal N}\times \CP^1$ by the twistor
distribution (Lax system) ${L_{Ai}}$. This clearly has a projection
$q :{\cal N} \times \CP^1 \mapsto {\cal PT}$ and we have a double
fibration
$$
\begin{array}{rcccl}
&&{\cal N}\times\CP^1&&\\
&p\swarrow&&\searrow q&\\
{\cal N}&&&&{\cal PT}
\end{array}
$$

Since the twistor distribution is tangent to the fibres of ${\cal
  N}\times \CP^1\mapsto \CP^1$, twistor space inherits the projection
$\mu:{\cal PT}\mapsto \CP^1$.  The twistor space for the hierarchy is
three-dimensional as for the ordinary hyper-K\"ahler equations, but has
a different topology.  We have

\begin{lemma} The holomorphic curves $q(\CP^1_x)$ where
$\CP^1_x=p^{-1}x$, $x\in {\cal N}$, have normal bundle $N={\cal
  O}(n)\oplus{\cal O}(n)$.  
\end{lemma}

\noindent
{\bf Proof.} To see this, note that $N$ can be identified with the
quotient $p^*(T_x{\cal N})/\{ \mathrm{span }L_{Ai}\}$, $i=1,\ldots,n$.
In their homogeneous form the operators $L_{Ai}$ have weight 1, so the
distribution spanned by them is isomorphic to the bundle
$\C^{2n}\otimes{\cal O}(-1)$.  The definition of the normal bundle as
a quotient gives
$$
0\mapsto \C^{2n}\otimes{\cal O}(-1) \mapsto \C^{2n+2}  \mapsto N\mapsto 0
$$
and we see, by taking determinants that the image is ${\cal
  O}(n+a)\oplus {\cal O}(n-a)$ for some $a$.  We see that $a=0$ as the
last map, in the spinor notation introduced at the end of the last
section, is given explicitly by $V^{AA'_1\ldots A'_n}\mapsto
V^{AA'_1\ldots A'_n}\pi_{A'_1} \ldots \pi_{A'_n}$ clearly projecting
onto ${\cal O}(n)\oplus {\cal O}(n)$.  \hfill \koniec

\medskip 

A final structure that ${\cal PT}$ possesses is a skew form $\Sigma$
taking values in ${\cal O}(2n)$ on the fibres of the projection $\mu$.
This arises from the fact that the vector fields of the distribution
preserve the coordinate volume form $\nu$ on ${\cal N}$ in the given
coordinates system.  Furthermore, the Lax system commutes exactly
$[L_{aI},L_{Bj}]=0$ so that
$$
\Sigma=\nu(\cdot,\cdot,L_{01},\ldots, L_{0n}, L_{11},\ldots ,L_{1n})
$$
descends to the fibres of ${\cal PT}\mapsto \CP^1$ and clearly has
weight $2n$ as each of the $L_{Ai}$ has weight one.

Thus we see that, given a solution to the hyperk{\"a}hler hierarchy in
the form of a commuting Lax system, we can produce a twistor space
with the above structures.  Now we shall prove the main result of this
section and demonstrate that, given ${\cal PT}$, with the above
structures, we can construct ${\cal N}$ (as the moduli space of
rational curves in ${\cal PT}$) which is naturally equipped with a
function $\Th$ satisfying (\ref{hier2}) and with the Lax distribution
(\ref{system}).

\begin{prop}
\label{suptw}
Let $\cal PT$ be a 3 dimensional complex manifold with the following 
structures
\begin{itemize}
\item[1)] a projection $\mu :{\cal PT}\longrightarrow \CP^1$,
%\item[2)] Euler homogeneity operator $
%\Upsilon=\om^A{\p}/{\p\om^A}+\pi_{A'}{\p}/{\p\pi_{A'}}$ on $\cal T$, 
\item[2)] a section $s:\CP^1\mapsto {\cal PT}$ of $\mu$
% homogeneous of degree $n$ with respect to $\Upsilon$, each 
  with normal bundle ${\cal O}(n)\oplus{\cal O}(n)$,
\item[3)] a non-degenerate 2-form $\Sm$
%=d\om^A\wedge d\om_A$ 
  on the fibres of $\mu$, with values in the pullback from $\CP^1$ of
  ${\cal O}(2n)$.
\end{itemize}
Let ${\cal N}$ be the moduli space of sections that are
deformations of the section $s$ given in 
{\em(2)}. Then ${\cal N}$ is $2n+2 $ dimensional and
\begin{itemize}
\item[a)] There exists coordinates, $x^{Ai}$, $A=0,1$, and $i=0,\ldots
, n$ and a function $\Th:{\cal N}\longrightarrow \C$ on ${\cal N}$
such that equation {\em(\ref{hier2})} is satisfied.
\item[b)] The moduli space ${\cal N}$ of sections is equipped with 
\begin{itemize}
\item a factorisation of the tangent bundle $T{\cal N}=S^A\otimes
  \odot^n S^{A'}$, 
\item a $2n$-dimensional distribution on the `spin bundle' $D\subset
  T({\cal N}\times \CP^1)$ that is tangent to the fibres of {\bf r}
  over $\CP^1$ and, as a bundle on ${\cal N}\times\CP^1$ has an
  identification with ${\cal O}(-1)\otimes S_{AA'...A_{n-1}'}$ so that
  the linear system can be written as in equation {\em(\ref{system})}.
\end{itemize}
\end{itemize}
This correspondence is stable under small perturbations of the complex
structure on ${\cal PT}$ preserving (1) and (3).
\end{prop}

\noindent
{\bf Proof:}  The first claim, that ${\cal N}$ has dimension $2n+2$
follows from Kodaira theory as $\dim H^0(\CP^1, N)=2n+2$  and $\dim
H^1(\CP^1, N)=\dim H^1(\CP^1, \mathrm{End} N)=0$.

\begin{figure}
\caption{ Double fibration. }
\label{dopic}
\begin{center}
\includegraphics[width=11cm,height=8cm,angle=0]{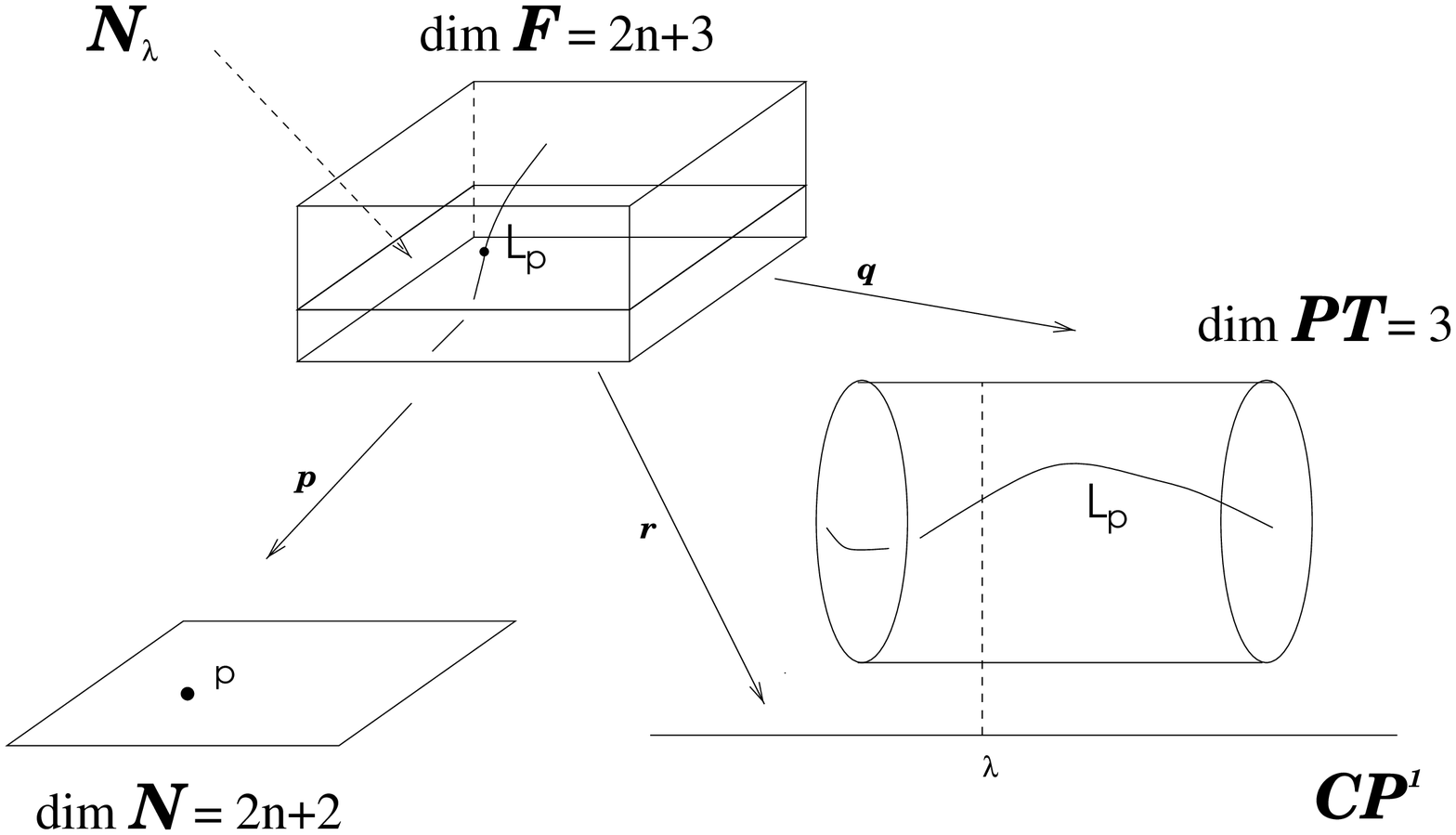}
\end{center}
\end{figure}

\smallskip

\noindent
{\bf Proof of (a):} we first start by defining homogeneous coordinates on
${\cal PT}$.  These are coordinates on ${\cal T}$, the total space of
the pullback from $\CP^1$ of the tautological line bundle ${\cal
O}(-1)$.  Let $\pi_{A'}$ be homogeneous coordinates on $\CP^1$ pulled
back to $\cal T$ and let $\omega^A$ be local coordinates on $\cal T$
chosen on a neighbourhood of $\mu^{-1}\{\pi_{0'}=0\}$ that are
homogeneous of degree $n$ and canonical so that
$\Sigma=\varepsilon_{AB} \d\omega^A \wedge \d\omega^B$.  We also use
$\lambda=\pi_{0'}/\pi_{1'}$ as an affine coordinate on $\CP^1$. Let
$L_p$ be the line in $\cal PT$ that corresponds to $p\in {\cal N}$ and
let $Z\in\cal PT$ lie on $L_p$.  We denote by $\cal F$ the
correspondence space ${\cal PT}\times {\cal N}|_{Z\in L_p}= {\cal
N}\times\CP^1$. (See figure 1 for the double fibration picture.)

Pull back the twistor coordinates to $\cal F$ and define 
$2(n+1)$ coordinates on ${\cal N}$ by
\[
x^{A(A_1'A_2'...A_n')}:=\frac{\p^n \om^A}{\p
\pi_{A_1'}\p\pi_{A_2'}...\p\pi_{A_n'}}{\Big |}_ {\pi_{A'}=o_{A'}},
\]
where the derivative is along the fibres of $\cal F$ over $\cal N$.
This can alternatively be expressed in affine coordinates on $\CP^1$
by expanding the coordinates $\omega^A$ pulled back to $\cal F$ in
powers of $\lambda=\pi_{0'}/\pi_{1'}$:
\be
\label{expan}
\om^A = (\pi_{1'})^n\left(
  \sum_{i=0}^{n}x^{Ai}\lambda^{n-i}+\lambda^{n+1}
  \sum_{i=0}^{\infty}s^{A}_i\lambda^{i} \right), 
\ee 
where the $s^A_i$ are functions of $ x^{AA_1'...A_n'}$ and will be
useful later.

The symplectic 2-form $\Sm$ on the fibres of $\mu$, when pulled back
to the spin bundle, has expansion in powers of $\lambda$ that
truncates at order $2n+1$ by globality and homogeneity, so that
\[
\Sm=\d_h\om_A\wedge \d_h\om^A=\pi_{A_1'}...\pi_{A_n'}\pi_{B_1'}...\pi_{B_n'}
\Sm^{A_1'...A_n'B_1'...B_n'}
\]
for some symmetric spinor indexed 2-form $\Sm^{A_1'...A_n'B_1'...B_n'}$.  We have
\be
\label{algformy}
\Sm(\lambda)\wedge\Sm(\lambda)=0,\;\;\;\;\;\d_h\Sm(\lambda)=0.  
\ee 
where in the exterior derivative $\d_h$,  $\lambda$ is understood to
be held constant.

If we express the forms in terms of the $x^{Ai}$ and the $s^A_i$, the
closure condition is satisfied identically, whereas the truncation
condition will give rise to equations on the $s^A_i$
allowing one to express them in terms of a function
$\Th(x^{AA'...A_n'})$ and to field equations on $\Th$ as follows. 

%Set
%\[
%\Sm^i={2n \choose i}\Sm^{A_1'...A_n'B_1'...B_n'} o_{A_1'}...o_{A_i'}
%\iota_{A_{i+1}'}...\iota_{B_{2n}'}(-1)^{2n-i},
%\]
%so that 
%$
%\Sm=\sum_{i=0}^{2n}\lambda^{2n-i}\Sm^i
%$
%where
%\begin{eqnarray*}
%\Sm^j & =&  \sum_{i=j-n}^{n}\varepsilon_{AB}dx^{Ai}
%\wedge \d x^{Bj-i}\;\;   \mbox{ for $j\geq n $} \\
%     & =&    \sum_{i=0}^{j}\varepsilon_{AB}\d x^{Ai}
%\wedge \d x^{Bj-i}+
%  \sum_{i=j+1}^{n}\varepsilon_{AB}\d x^{Ai}
%\wedge \d s^{B}_{i+j-1}\;  \mbox{ for $j< n$ }\nonumber.
%\end{eqnarray*}

%The algebraic condition in (\ref{algformy}) implies
%\be
%\label{roow}
%\sum_{i=\mbox{max }\{0,k-2n\}}^{\mbox{min} \{2n,k\}}
%\Sm^i\wedge\Sm^{k-i}=0\;\;\;\;\;\;\;\;\;\;\;\;\;      
%\mbox{for}\;0\leq k\leq 4n.
%\ee

To deduce the existence of $\Th(x^{AA_1'...A_n'})$ observe that the
vanishing of the coefficient of $\lambda^{2n+1}$ in $\d\om^A\wedge
\d\om_A$ gives
\[
\sum_{i=0}^{n}\d s_{Ai}\wedge \d x^{Ai}= \d \sum_{i=0}^{n}
s_{Ai} \d x^{Ai}= 0 \qquad \Longrightarrow \qquad
s_{Ai}=\frac{\p \Th}{\p x^{Ai}}.
\]

The equations of the hierarchy arise from
the vanishing of the coefficient of $\lambda^{2n+2}$
$$
\sum_{i=0}^n \d x^{Ai}\wedge \d s_A^{i+1} + \d s^{A0}\wedge \d
s_A^{0}=0 \, .
$$
This leads to the equations (\ref{hier2}) on $\Theta$ for $i,j\leq
n-1$,
$$
\frac{\p^2 \Theta}{\p x^{Ai+1}\p x^{Bj}} - \frac{\p^2 \Theta}{\p
  x^{Ai}\p x^{Bj+1}} + \varepsilon^{CD}\frac{\p^2 \Theta}{\p x^{C0}\p
  x^{Ai}}\frac{\p^2 \Theta}{\p x^{D0}\p x^{Bj}}=0
$$
and further equations that determine $s^{An+1}$.

% which is a particular case of (\ref{roow})
%\begin{eqnarray*}
%0 &=& 2\sum_{j=1}^{n} \sum_{i=1}^{n}\varepsilon_{AB}\varepsilon_{CD}
%\d x^{C}\wedge \d x^{D}\wedge \d x^{Ai}\wedge \d x^{Bj}
%\frac{\p^2\Th}{\p x^{Ai-1}\p x^{Bj}}\\
%&+&\sum_{i=1}^{n} \sum_{k=0}^{n}\sum_{j=1}^{n} \sum_{l=0}^{n}
%\d x^{Ai}\wedge \d x^{Ck}\wedge \d x^{Bj}\wedge \d x^{Dl}
%\frac{\p^2\Th}{\p x^{Ai-1}\p x^{Ck}}\frac{\p^2\Th}{\p x^{Bj-1}\p x^{Dl}}.
%\end{eqnarray*}
%Equate the coefficient of 
%$\varepsilon_{AB}\d x_C\wedge \d x^C\wedge \d x^{Ai}\wedge \d x^{Bj}$ to
%$0$ (i.e.\ the term with $k=l=0$). This yields (\ref{hier2})
%\be
%\p_{Ai}\p_{Bj-1}\Th -\p_{Bj}\p_{Ai-1}\Th +\{ \p_{Ai-1}\Th, \p_{Bj-1}\Th
%\}_{yx}=0,
%\ee
%and all other equations in (\ref{roow}) are satisfied trivially.

\medskip

\noindent
{\bf Proof of b).} 
The isomorphism $T{\cal N}=S^A\otimes \odot^n S^{A'}$ follows simply
from the structure of the normal bundle. From Kodaira theory,
since the appropriate obstruction groups vanish, we have
\begin{equation}\label{solder}
T_x{\cal N}=\Gamma(\CP^1_x, N_x)=S_x^A\otimes \odot^n S^{A'} 
\end{equation}
where $N_x$ is the normal bundle to the rational curve $\CP^1_x$ in
$\cal PT$ corresponding to the point $x\in \cal N$.  The bundle $S^A$
on space-time is the Ward transform of ${\cal O}(-n)\otimes T_V{\cal
  PT}$ where the subscript $V$ denotes the sub-bundle of the tangent
bundle consisting of vectors up the fibres of $\mu$, the 
projection to $\CP^1$, so that
$S^A_x=\Gamma(\CP^1_x,{\cal O}(-n)\otimes T_V{\cal PT})$.  The bundle
$S^{A'}=\Gamma(\CP^1,{\cal O}(1))$ is canonically trivial.  

Let
$\nabla_{AA'_1\cdots A'_n}=\nabla_{A(A'_1\cdots A'_n)}$ be the indexed
vector field that establishes the isomorphism (\ref{solder}) and let
$e^{AA'_1\cdots A'_n}=e^{A(A'_1\cdots A'_n)}\in \Omega^1\otimes S^A\otimes
\odot^n S^{A'}$ be the dual (inverse) map. 

We now wish to derive the form of the linear system, equations
(\ref{system}).  For each fixed $\pi_{A'}=(\lambda, 1)\in \CP^1$ we
have a copy of a space-time ${\cal N}_{\lambda}$. The horizontal
(i.e.\ holding $\lambda$ constant) subspace of $T_{(x,\lambda)}({\cal
  N}\times \CP^1)$ is spanned by $\nabla_{A(A'...A_n')}$. An element
of the normal bundle to the corresponding line $\CP^1_x$ consists of a
a horizontal tangent vector at $(x, \lambda)$ modulo the twistor
distribution. Therefore we have the sequence of sheaves over $\CP^1$
%Lax distribution is $\pi^{A'}\nabla_{AA'A_2'...A_n'}\in 
%S_{AA'...A_n'}\otimes{\cal O}(1)\otimes D$
\[
0\longrightarrow D_x \longrightarrow T_x{\cal N}\stackrel{e^A} \longrightarrow
S^A\otimes{\cal O}(n)\longrightarrow 0\, ,
\]
where $D_x$ is the twistor distribution at $x$ and the map $T_x{\cal
  N} \longrightarrow S^A\otimes{\cal O}(n)$ is given by the
contraction of elements of $T_x{\cal N}$ with
$e^A:= e^{AA_1'...A_n'}\pi_{A_1'}...\pi_{A_n'}$ since $e^A$ annihilates
all $L_{Bi}$s in $D$.  Consider the dual sequence
%\[
%0\longrightarrow {\cal O}_A(-n) \longrightarrow T^{*}{\cal N} \longrightarrow
%D^*\longrightarrow 0
%\]
tensored with ${\cal O}(-1)$ to obtain
\be
\label{short}
0\longrightarrow {\cal O}_A(-n-1) \longrightarrow T^{*}_x{\cal N}(-1)
 \longrightarrow D^*_x(-1)\longrightarrow 0.
\ee

From here we would like to extract the Lax
distribution  $$L_{AA_2'...A_n'}=\pi^{A_1'}D_{AA_1'A_2'...A_n'}
\in S_{AA_2'...A_n'}\otimes{\cal O}(1)\otimes D.$$
This can be achieved by globalising (\ref{short}) in $\pi^{A'}$ .
The corresponding long exact sequence of cohomology groups yields
\[
0\longrightarrow \Gamma({\cal O}_A(-n-1)) \longrightarrow
\Gamma(T^{*}{\cal N}(-1))
\longrightarrow \Gamma(D^*(-1))\stackrel{\delta}{\longrightarrow} 
H^1({\cal O}_A(-n-1))
\]
\[ \longrightarrow H^1(T^{*}{\cal N}(-1)) \longrightarrow
... 
\]
which (because $T^*\cal N$ is a trivial bundle so that ${\cal
  O}(-1)\otimes T^*\cal N$ has no sections or cohomology) reduces to
\[
0\longrightarrow \Gamma(D^*(-1))\stackrel{\delta}{\longrightarrow}
H^1({\cal O}_A(-n-1)) \longrightarrow 0.  
\] 
From Serre duality 
we conclude, since $D$ has rank $2n$, that the connecting map $\delta$
is an isomorphism $\delta :\Gamma(D^*(-1))\longrightarrow S_{AA_2'...A_{n}'}$.
Therefore
\be
\delta \in \Gamma(D\otimes {\cal O}(1)\otimes S_{AA_2'...A_{n}'})
\ee
is a canonically defined object annihilating $\om^A$ given by
(\ref{expan}). 

In index notation we can put
\[
\delta =L_{AA_2'...A_n'}
=\pi^{A_1'}D_{AA_1'A_2'...A_n'},
\]
where $L_{AA_2'...A_n'}=L_{A(A_2'...A_n')}$, the second identity
follows from the globality of $L_{AA_2'...A_n'}$ and the
$D_{AA_1'A_2'...A_n'}$ are vector fields on $\cal N$ lifted to ${\cal
  N}\times \CP^1$ using the product structure. 

It follows from $L_{AA_2'...A_n'}\om^B=0$ that if $\pi^{A'}=o^{A'}$
then $D_{A0'A_2'...A_n'}x^{Bn}=0$ so
\[
D_{A0'A_2'...A_n'}=A_{A0'A_2'...A_n'}^{BB_2'...B_n'}
\frac{\p}{\p x^{B0'B_2'...B_n'}}\, ,
\]
for some matrix $A_{A0'A_2'...A_n'}^{BB_2'...B_n'}$.  This matrix must
be invertible by dimension counting.  By multiplying
$L_{AA_2'...A_n'}$ by the inverse of this matrix, we find we can put
\[
A_{A0'A_2'...A_n'}^{BB_2'...B_n'}=
\varepsilon_A^B\varepsilon_{A_2'}^{B_2'}...\varepsilon_{A_n'}^{B_n'}.
\]
Therefore we can take $
L_{AA_2'...A_n'}=\p_{A0'A_2'...A_n'}-\lambda D_{A1'A_2'...A_n'}$.
Equating the $(n-i+1)$th and $(n+1)$th powers of $\lambda$ in
$L_{Ai}\om^B=0$ to zero yields
\[
D_{A1'A_2'...A_n'}=\p_{A1'A_2'...A_n'}+[\p_{A0'A_2'...A_n'}, V]
\]
where $V=\varepsilon_{AB}\p\Th/\p x_{A0} \p/\p x_{B0}$. So
finally $L_{AA_2'...A_n'}$ is of the form 
$L_{Ai}=\p_{Ai}-\lambda(\p_{Ai+1}+[\p_{Ai}, V])$.
\koniec

\subsection{Geometric structures}
\label{geomstr}
If one considers ${\cal N}={\cal M}\times \X$ as being foliated by
four dimensional slices $t^{Ai}=const$ then structures (1)--(3) on
${\cal PT}$ can be used to define anti-self-dual vacuum metrics on the
leaves of the foliation.  Consider $\Th(x^{AA'}, {\bf t})$ where ${\bf
t}=\{ t^{Ai}, i=2...n\}$. For each fixed ${\bf t}$ the function $\Th$
satisfies the second heavenly equation. The ASD metric on a
corresponding four-dimensional slice ${\cal N}_{{\bf t}={\bf t_0}}$ is
given by
\[
\d s^2=2\varepsilon_{AB}\d x^{A1'}\d x^{B0'}+
2\frac{\p^2 \Th}{\p x^{A0'}\p x^{B0'}}\d x^{A1'}\d x^{B1'}.
\]
This metric can be determined from the structure
of the ${\cal O}(n)\oplus{\cal O}(n)$ twistor space as follows.

Fix the first $2n-2$ parameters in the expansion (\ref{expan}) so the
normal vector $W=W^A\p/\p\om^A$ is given by
\[
W^A=\delta \om^A=\lambda ^{n-1}W^{A1'}+\lambda ^{n}W^{A0'}+
\lambda^{n+1}\frac{\p \delta\Th}{\p x_{A}^{0'}}+...
\]
where $\delta\Th=W^{AA'}\p\Th/\p x^{AA'}$. The metric is
\be
\label{cnstructure}
g(U,W)=\varepsilon_{AB}\varepsilon_{A'B'}U^{AA'}W^{BB'}
\ee
where $\varepsilon_{A'B'}$ is a fixed element of $\Lambda^2S^{A'}$ and 
$\varepsilon_{AB}\in \Lambda^2S^{A}$ is determined by $\Sigma$; recall
that $S^A_x=\Gamma (L_x, {\cal O(-n)}\otimes T_V{\cal
PT})$.  Thus if $u^A, v^A\in S_x^A$, then define
$\varepsilon_{AB}u^Av^B=\Sigma(u,v)$ where $u, v$ are the
corresponding weighted vertical vector fields on ${\cal PT}$. 

%\Here $\alpha_{A'}$ and $\beta_{A'}$ are zeros of $U$ and $W$.
%The last formula follows also from  (\ref{parakon}) if one puts
%$$W^{AA_1'...A_n'}=W^{A(A_1'}o^{A_2'}...o^{A_n')}$$ for $W$ tangent to 
%$t^{Ai}=const$.

%Note that it is sufficient to consider the slice ${\bf t}=0$.
%This is because an appropriate (canonical) coordinate transformation
%of $\cal PT$, $\om^A\rightarrow \hat {\om}^A(\om^B, \lambda)$  induces the
%transformation of parameters $
%\{{\bf t}={\bf t_0} \} \rightarrow \{{\bf \hat t}=0\}$.

For $n$ odd $T\cal N$ is equipped with a metric with
holonomy $SL(2,\C)$.  For $n$ even, $T\cal N$ is endowed with a skew
form. They are both given by
\be
\label{parakon}
G(U,
W)=\varepsilon_{AB}\varepsilon_{A_1'B_1'}...\varepsilon_{A_n'B_n'}
U^{AA_1'...A_n'}W^{BB_1'... B_n'}.  
\ee 
These are special examples of the paraconformal structures considered
by Bailey and Eastwood \cite{BE91}.

\subsection{Holomorphic deformations and $\O(2n)$ twistor functions} 
\label{infdef}
%From the discussion in the beginning of Subsection 
%\ref{ONONtwistor} if follows that the

We wish to consider holomorphic deformations of ${\cal PT}$ that
preserve conditions $(1-3)$ of Proposition \ref{suptw} which will
therefore correspond to perturbations of the hierarchy.

Let $\tom^A=G^A(\om^B,\pi_{A'},t)$ be the standard patching relation
for $\Pt$ and let $f^A\in S^A\otimes H^1(\Pt,\O(n))$ give the
infinitesimal deformation
\[
\tom^A=G^A+ t f^A +O(t^2).
\]
The globality of the symplectic structure 
$\d\tom_A\wedge \d\tom^A=\d\om_A\wedge \d\om^A$
implies $f^A=\varepsilon^{AB}\p f/\p \om^B$ 
where $f\in H^1(\Pt, \O(2n))$.

\smallskip

\noindent
{\bf Example:} if we deform from the flat model using  
$f=(\pi_{0'})^{4n}/\om^0\om^1$,
then the deformation equations
\[
\tom^0=\om^0+t\frac{(\pi_0)^{4n}}{\om^0(\om^1)^2}+O(t^2),\qquad
\tom^1=\om^1-t\frac{(\pi_0)^{4n}}{(\om^0)^2\om^1}+O(t^2).
\]
imply that $Q=\om^0\om^1=\tom^0\tom^1$ is a global twistor function
(up to $O(t^2)$) which persists to all orders as 
$\varepsilon^{AB}\p Q/\p\om^A\p f/\p \om^B=0$.  
The corresponding deformed paraconformal
structure admits a symmetry corresponding to the global vector field
$\varepsilon^{AB}\p Q/\p\om^A\p/\p \om^B$ on ${\cal PT}$.

\smallskip

To see how such `Hamiltonians' $f$ correspond to variations in the
paraconformal structure (or more simply $\Theta$), we form an indexed
element of $H^1(\Pt, \O(-1))$, and pull it back to 
${\cal N}\times\CP^1$ where it can be split uniquely:
\[
\pi_{A_2'}...\pi_{A_n'}\frac{\p^3f^{2n}}{\p\om^{A}\p\om^{B}\p\om^{C}}
=f_{ABCA_2'...A_n'}
= \tilde{\cal F}_{ABCA_2'...A_n'}-{\cal F}_{ABCA_2'...A_n'}.
\]
where
\[
{\cal F}_{ABCA_2'...A_n'}=\frac{1}{2\pi i} \oint_{\Gamma}
\frac{f_{ABCA_2'...A_n'}}{\rho_{A'}\pi^{A'}}\rho
\cdot \d \rho.
\]
This gives rise to a global field that is symmetric over its indices:
\[
C_{ABCDA_2'...A_n'D_2'...D_n'}=L_{DD_2'...D_n'}{\cal F}_{ABCA_2'...A_n'}
\]
which is given also directly by the integral 
\[
C_{ABCDA_2'...A_n'D_2'...D_n'}=
\frac{1}{2\pi i} \oint_{\Gamma}
\rho_{A_2'}...\rho_{A_n'}\rho_{D_2'}...\rho_{D_n'}
\frac{\p^4f^{2n}}{\p\om^{A}\p\om^{B}\p\om^{C}\p\om^{D}}
\rho\cdot\d\rho.
\]

To see how this corresponds to a variation of $\Theta$, we introduce a
chain of potentials.  Use the non-unique splitting $f^{2n}={\cal
  F}^{2n}-\widetilde{\cal F}^{2n}$ and define a global object of
degree $2n+1$ by
\[
L_{AA_2'...A_n'}{\cal F}^{2n}=\Sm_{AA_2'...A_n'B_1'...B_n'C_1'...C_n'D_1'}
\pi^{B_1'}...\pi^{B_n'}\pi^{D_1'}\pi^{C_1'}...\pi^{C_n'}.
\]
It is easy to see that
\[
\nabla^{AE_1'...E_n'}\Sm_{AA_2'...A_n'B_1'...B_n'C_1'...C_n'D_1'}=0,
\]
and $\Sm_{AA_2'...A_n'B_1'...B_n'C_1'...C_n'D_1'}$ is a potential
potentials, related to the field by
\[
C_{ABCDA_2'...A_n'D_2'...D_n'}=
\nabla^{D_1'}_{DD_2'...D_n'}
\nabla^{C_1'...C_n'}_C
\nabla^{B_1'...B_n'}_B\Sm_{AA_2'...A_n'B_1'...B_n'C_1'...C_n'D_1'}.
\]
%Use the gauge where all potentials with non zero number of primed
%indices vanish, when contracted with $o^{A'}$
The chain of potentials is
\begin{eqnarray*}
\ddot_{A_1'B_1'...B_n'C_1'...C_n'D_1'}
&=&o_{A_1'}o_{B_1'}...o_{B_n'}
o_{C_1'}...o_{C_n'}o_{D_1'}\ddot\\
\Sm_{AA_2'...A_n'B_1'...B_n'C_1'...C_n'D_1'}
&=&o_{B_1'}...o_{B_n'}o_{C_1'}...o_{C_n'}o_{D_1'}
\nabla_{A0'A_2'...A_n'}\ddot\\
H_{ABA_2'...A_n'B_1'...B_n'D_{1'}}&=&o_{B_1'}...o_{B_n'}o_{D_1'}\nabla_{B0'}\nabla_{A0'A_2'...A_n'}\ddot\\
\Gamma_{ABCA_2'...A_n'D_1'}
&=&o_{D_1'}\nabla_{C0'}\nabla_{B0'}\nabla_{A0'A_2'...A_n'}\ddot\\
C_{ABCDA_2'...A_n'D_2'...D_n'}&=&
\nabla_{C0'}\nabla_{B0'}\nabla_{A0'A_2'...A_n'}\nabla_{D0'D_2'...D_n'}\ddot.
\end{eqnarray*}
This can be compared with the corresponding chain for $n=1$ \cite{KLNT81}.
%%%%%%%%%%%%%%%%%%%%%%%%%%%%%%%%%%%%%%%%%%%%%%%%%%%%%%%%%%%%%%%%%%%%%%%%%%%%%%

\section{Hamiltonian and Lagrangian  formalisms}
\label{hformalism}

In this Section we shall investigate the Lagrangian and Hamiltonian
formulations of the hyper-K\"ahler equations in their `heavenly' forms.
The symplectic form on the space of solutions to heavenly
equations will be derived, and proven to be compatible with a
recursion operator.

Both the first and second heavenly equations admit Lagrangian
formulations, and these can be used to derive symplectic structures on
the solution spaces, which we denote by $\cal S$.  Here, rather than
consider the equations as a real system of elliptic or ultra-hyperbolic
equations, we complexify and consider the equations locally as
evolving initial data from a 3-dimensional hyper-surface and it is this
space of initial data that leads to local solutions on a neighbourhood
of such a hyper-surface that is denoted by $\cal S$ and is endowed with a
(conserved) symplectic form.

For the first equation we have the Lagrangian density
\be
\label{lagr1}    
{\cal L}_{\Om}=\Om \Big(\nu - \frac{1}{3}(
\p\tilde{\p}
\Om)^2\Big)=\Big(\Om-\frac{1}{3}\Om\{\Om_{\zt}, \Om_{\wt} \}_{wz}\Big)\nu
\ee
and for the second equation 
\begin{eqnarray}
\label{lagr2}
{\cal L}_{\Th}&=&\Big(\frac{2}{3}
\Th(\dd\Th)^2-
\frac{1}{2}
(\p\Th)\wedge(\dd \Th)\Big)\wedge e^{A0'}\wedge e_A^{0'}\nonumber\\
&=&\Big(\frac{1}{3}\Th \{\Th_x,\;
\Th_y\}_{xy}-\frac{1}{2}(\Th_x\Th_w+\Th_y\Th_z)\Big)
\nu.
\end{eqnarray}
Note that $e^{A0'}\wedge e_A^{0'}$ can be replaced by
$\d x\wedge \d y$ in the second Lagrangian as it is multiplied by
$\d w\wedge \d z$. 

If the field equations are assumed, the variation of these Lagrangians
will yield only a boundary term.  Starting with the first equation,
this defines a potential one-form $P$ on the solution space $\cal S$
and hence a symplectic structure ${\bf \Om} =\d P$ on $\cal S$.
Starting with the second we find a symplectic structure with the same
expression on perturbations $\delta \Th$ as we had for $\delta\Omega$.
However, since their relation to perturbations of the hyper-K\"ahler
structure are different, they define different symplectic structures
on $\cal S$.  These are related by the recursion operator since we have
$R^2\delta \Omega = \delta \Theta$ from above.  In order to see that
these structures yield the usual bi-Hamiltonian framework, we will
need to show that these symplectic structures are compatible with the
recursion operator in the sense that ${\bf \Om} (R\phi,\phi') ={\bf
\Om} (\phi,R\phi')$. 

We shall demonstrate this using the first heavenly formulation which
is easier as one can use identities from K\"ahler geometry.  (The
derivation of the symplectic structure from the second Lagrangian will
be done in coordinates, since the useful relation between the Hodge
star and the K\"ahler structure is missing in this case.)

\begin{prop}
\label{o_symplek}
The symplectic form on the space of solutions $\cal S$ derived
from the boundary term in the variational principle for the first
Lagrangian is 
\be
\label{sympstruc}
{\bf \Om}(\delta_1\Om, \delta_2\Om) =\frac{2}{3}\int_{\delta M} \delta _1\Om \ast \d(\delta _2 \Om) -
\delta _2 \Om \ast \d(\delta _1\Om ).
\ee
\end{prop}
{\bf Proof.} Varying (\ref{lagr1}) we obtain
\[
\delta L =\delta \Om (\nu - \frac{1}{3}(\p \pt \Om)^2) -
\frac{2}{3}\Om \p \pt \Om\wedge  \p \pt \delta \Om
=\frac{2}{3}\p \pt \Om\wedge (\delta \Om\p \pt \Om -\Om \p \pt \delta \Om).
\]
%\[
%= dA(\delta \Om) +(Euler-Lagrange)\delta \Om.
%\]
We use the identities $\d(\p - \pt ) = 2\pt \p,\;\;\;
\omega \wedge J_1\d = \p \pt \Om \wedge (\p - \pt ) = \ast \d$ and the
field equation to obtain
\begin{eqnarray*}
 \delta L &=& -\frac{1}{3}\p \pt \Om \wedge (\delta \Om \d(\p - \pt )\Om
- \Om \d(\p - \pt ) \delta \Om )\\
%&=& -\frac{1}{3}\p \pt \Om \wedge( \d(\delta \Om (\p - \pt)\Om 
%-\Om (\p - \pt)\delta \Om) 
% -\d(\delta \Om )\wedge(\p -\pt )\Om + \d\Om \wedge (\p -\pt )\delta \Om )\\
%&=& -\frac{1}{3}\d(\delta \Om \ast \d\Om - \Om \ast \d(\delta \Om ))
% - \frac{1}{3}\p \pt \Om \wedge (-\ast \ast \d(\delta \Om )
%\wedge (\p -\pt )\Om
%+ \ast \ast \d\Om\wedge (\p - \pt )\delta \Om )\\
&=& \frac{1}{3} \d A(\dom ) - \frac{1}{3}\p \pt \Om (-\ast \p \pt \Om (\p -\pt )
\dom (\p - \pt ) \Om 
 + \ast \p \pt \Om (\p -\pt ) \Om (\p -\pt )\dom )\\
&=& \frac{1}{3} \d A(\dom ) \;\;\;\; \mbox{where}\;\;
A(\dom )=\Om \ast \d\delta \Om -\delta \Om \ast \d\Om. 
\end{eqnarray*}
Define the one form on $\cal S$
\[
P =\int_{\delta M} A(\dom ).
\]
The symplectic structure ${\bf{\Om}}$ is the (functional) exterior
derivative of $P$
\begin{eqnarray*}
{\bf{\Om}}(\delta _1 \Om ,\delta _2 \Om ) &=&\delta _1 (P(\delta _2 \Om )) -
\delta _2 (P(\delta _1 \Om )) -P([\delta _1 \Om , \delta _2 \Om ])\\
&=&\frac{2}{3}\int_{\delta M} \delta _1\Om \ast \d(\delta _2 \Om) -
\delta _2 \Om \ast \d(\delta _1\Om ). \qquad\qquad\qquad \Box
\end{eqnarray*}

%By a formal application of Stokes' theorem
%\[
%{\bf \Om}(\delta_1\Om, \delta_2\Om)
%=\frac{2}{3}\int_{\delta {\cal M}}{e^{1'}}_B\wedge e^{B0'}\wedge
%(e^{A0'} \delta _2\Om\nabla_{A0'}\delta_1\Om - 
%e^{A1'} \delta _1\Om\nabla_{A1'}\delta_2\Om)
%\] 

Thus ${\bf \Om}$ coincides with the symplectic form on the solution
space to the wave equation on the ASD vacuum background.
 
The existence of the recursion operator allows the
construction of an infinite sequence of symplectic structures.
The key property we need is the following
\begin{prop}
\label{o_bih}
Let $\phi ,\;\phi' \in W_g$ and let ${\bf \Om}$ be given by
{\em(\ref{sympstruc})}. Then
\be
\label{olver}
{\bf \Om}(R\phi ,\;\phi')={\bf \Om}(\phi,\;R\phi').
\ee
\end{prop}
We first prove a technical lemma:
\begin{lemma}
The following identities hold
\begin{eqnarray}
\label{techlemma}
\om \wedge \p \phi &=&-\a \wedge \pt R\phi,\;\;\;
\om \wedge \dd \phi =\at \wedge \dd R\phi,\\
\om \wedge \dd R \phi &=&-\a \wedge \p_0 \phi,\;\;\;
\om \wedge \pt R \phi =\at \wedge \p \phi.\nonumber
\end{eqnarray}
\end{lemma}
{\bf Proof.} From the definitions of $\Sm^{A'B'}$ and $\p^{B'}_{A'}$ it follows
that 
\be
\Sm^{A'B'}\wedge \p^{C'}_{D'}=\Sm^{A'[B'}_{}\wedge \p^{C']}_{D'}
\ee
(recall that $\p^{B'}_{A'}=e^{AB'}\otimes\p_{AA'}$) which yields
\[
\om \wedge \pt =\at \wedge \dd ,\;\; \om \wedge \p = -\a \wedge \p_0,
\]
\[
\om \wedge \p_0 = \at \wedge \p ,\;\; \om \wedge \dd =-\a \wedge
\dd,\;\;\a \wedge \p =\at \wedge \pt=0.
\]
Multiplying (\ref{def}) by combinations of  spin co-frame  we get an
equivalent definition of the recursion operator
\be
\label{def3}
\p^{A'}_{1'}\phi=\p^{A'}_{0'}R\phi
\ee
which is equivalent to $\p \phi =\dd R \phi$ or $\do \phi= \pt R
\phi$. These formulae give the desired result.\koniec

\noindent
{\bf Proof of Proposition \ref{o_bih}.} 
The proof uses a (formal) application of Stokes' theorem:
\begin{eqnarray}
{\bf \Om}(\phi ,\;\phi')&=&\int_{\delta M} \phi \ast \d\phi' -
\phi' \ast \d\phi\nonumber\\
&=&\int_{\delta M} \om \wedge (\phi \p \phi' -\phi \pt \phi' -
\phi' \p \phi +\phi' \pt \phi )
=\int_{\delta M} \om \wedge (\phi \d \phi' +\phi' \d \phi' -
2\phi \pt \phi' -2\phi' \p \phi )\nonumber\\
\label{sto}
&=& -2\int_{\delta M} \om \wedge (\phi \pt \phi' +\phi' \p \phi)
=2\int_{\delta M}\om \wedge (\phi' \pt \phi +\phi \p
\phi')\nonumber.
\end{eqnarray}
From (\ref{sto}) and from (\ref{techlemma}) we have
\[
{\bf \Om}(\phi,\;R\phi') = -\int_{\delta M} \om \wedge 
(\phi \pt R \phi' +R\phi' \p \phi)
=-\int_{\delta M} \phi \p \phi' \wedge \at +\int_{\delta M} R\phi' \pt
R\phi \wedge \a
\]
and analogously
\[
{\bf \Om}(R\phi,\;\phi') =
 \int_{\delta M} \phi' \p \phi \wedge \at -\int_{\delta M} R\phi \pt
R\phi' \wedge \a.
\]
Equality (\ref{olver}) is achieved by subtracting the integral
of $\d (\phi \phi')\wedge \at -\d (R\phi R\phi')\wedge \a$ and
applying  Stokes' theorem. \koniec
\smallskip

This property guarantees that the bilinear forms 
\be
\label{struktury}
{\bf \Om}^k(\phi ,\;\phi')\equiv{\bf \Om}(R^k\phi ,\;\phi') 
\ee 
are skew.  Furthermore they are symplectic and lead to the
bi-Hamiltonian formulation.  In this context formula (\ref{olver})
and the closure condition for ${\bf \Om}^k$ are 
an algebraic consequence of the fact that $R$ comes from two Poisson
structures. Using the theory of bi-Hamiltonian systems one can now go
on to prove that the flows constructed by application of $R$ to some
standard flow commute. 

%Now we would like to use a sequence of symplectic structures to 

To develop the bi-Hamiltonian theory, we would like to write the
heavenly equations in Hamiltonian form. However the Legendre transform
becomes singular for the coordinate flows associated to the
coordinates we have chosen since they are, at least in the Minkowski
space limit, null coordinates.  One possibility is to develop a
Hamiltonian formalism based on such null hyper-surfaces. We shall
adopt a different approach and reformulate the second heavenly
equation as a first order system.  

%This will allow techniques known from the theory of bi-Hamiltonian
%systems. 
 
Define $\phi:=-\Th_x$ and formally rewrite the second heavenly equation
(\ref{secondeq}) as
\be
\label{firstorder}
\p_w\phi={\cal R}(\p_y\phi)\;\;\;\;\;\;\;\;\;
\mbox{where}\;\;\;\;\;
{\cal R}= (\p_z+\{\phi, ...\}_{yx})\circ {\p_x}^{-1}=
\nabla_{11'}\circ{\nabla_{10'}}^{-1}.
\ee
It is therefore a conjugated operator ${\cal R}$ (defined by 
(\ref{conjugate})), acting on
solutions to the zero-rest-mass equations, and plays the role of 
the recursion operator. Flows of the sub-hierarchy
$[L_1, L_{0j}]=0$  are
\[
\p_{t_j}\phi={\cal R}^j\p_y\phi
\]
and the Hamiltonian for the first nontrivial flow is
\[
H_1=\int\frac{\phi^2}{2} \d x \wedge \d y\wedge \d z.
\]
Higher Hamiltonians $H_n$ can in principle be constructed 
using the operator $R$.  However, we have not developed explicit
formulae for these $H_n$.  

\subsection{A local bi-Hamiltonian form for the hierarchy}
To end this section, we express the equations of the second heavenly
hierarchy (\ref{hier2}) in a compact form, and then write it as a
(formal) bi-Hamiltonian system on the spin bundle.  This will be a
rather different framework from that given above in that the
Hamiltonian structure will in effect be local to the $x^{A0}$ plane as
opposed to a field theoretic formulation---it is the gravitational
analogue of that given for the Bogomolny equations in \cite{MS92}
except that no symmetries are required here (in effect because ASD
gravity can be expressed as ASD Yang-Mills with two symmetries but
with gauge group the group of area preserving diffeomorphisms).  This
formulation is therefore presented merely as a curiousity.

Define the $j$th truncation of $\omega^A$ to be
\[
\om^A_j=-x^{A0}+ \sum_{m=1}^{j}\lambda^{m}\p^{Am-1}\Th
\, ,
\]
where $\p^{Ai}=\varepsilon^{AB}\p/\p x^{Bi}$.  (Note that this is
truncated at both ends, although the truncation at the lower end and
multiplication by a power of $\lambda$ is inessential.)
\begin{lemma}
The truncated heavenly hierarchy is equivalent to 
\be
\label{sato}
\p^{Bj} \;\om^A(\lambda)=
\{\om^A(\lambda), \lambda^{-j}\om^{Bj}(\lambda) \}_{yx}.
\ee
\end{lemma}
{\bf Proof.}
First  observe that one can sum the Lax system to obtain
\begin{eqnarray*}
-\sum_{i=0}^{j-1}\lambda^{i}L_{Ai}&=&\lambda^{j}\p_{Aj}+
\sum_{i=0}^{j-1}\lambda^{i+1} 
\varepsilon^{CD}\p_{C0}\p_{Ai}\Theta\p_{D0} 
-\p_{A0}\\
&=&\lambda^{j}\p_{Aj}+\{\omega_{Aj} \, , \cdot \}_{yx}
\end{eqnarray*}
where $\{f,\cdot\}_{yx}=\varepsilon^{CD}\p_{C0}f\p_{D0}$.

Thus, since $L_{Ai}\omega^A=0$, we have
$$
\p_{Bj}\omega^A= -\lambda^{-j}\{\omega^B_j\, , \omega^A\}
$$
which yields the desired answer. \koniec

%as a consequence of (\ref{hier2}),
%\[
%\sum_{m=0}^{j}\{\p_{Ai+j-m-1}\Th, \p_{Bm-1}\Th
%\}_{yx}=\p_{Bj}\p_{Ai-1}\Th.\]
%Therefore
%\[
%\p_{Bj}\sum_{i=0}^{\infty} \p_{Ai-1}\Th\lambda^{i+1}
%=\lambda^{-j-1} \sum_{a=0}^{\infty}\sum_{m=0}^{j}
%\{ \p_{Aa-1}\Th, \p_{Bm-1}\Th \} _{yx} \lambda^{a+m+2}.
%\]
%Fix $i$ and equate different powers of $\lambda$ on both sides.
%This  gives $a=i+j-m$ and finally (\ref{sato}).\koniec

\smallskip

%A functional derivative of any solution to Lax Pair 
%should be a Casimir function of the 
%Poisson pencil defined by 
%\be
%\{f, g \}_{\lambda}=\lambda\{f, g \}_i -\{f, g \}_{i-1}
%\ee

\noindent
For the remainder of this section, we shall fix the values of the 
spinor indices to be $A=0$ and $B=1$. 
Set 
\[
\p_j:=\p_{1j},\;\Psi:=\om^0,\;\mbox{and}\;\psi_j:=\om_{1j}.
\] 
Equation (\ref{sato}) takes the form
\[
\p_j\Psi(\lambda)=
\{ \Psi(\lambda), \lambda^{-j}\psi_j(\lambda) \} _{yx}
\]
which we rewrite as 
\be
\p_j\Psi={\cal D}\frac{\delta h_{j}}{\delta \Psi}.
\ee
Here ${\cal D}:=\{ \Psi(\lambda), ...\} _{yx}=
\sum_{i=0}^{\infty}{\cal D}_m\lambda^m$ is  $\lambda$-dependent 
Poisson structure, ${\cal D}_0=\p_x$ and 
${\cal D}_m=[\p_{m-1}, V]=D_{0m}-\p_{0m}$ for $m>0$. 

The Hamiltonians are
\[
h_j(\lambda)=\lambda^{-j}
\psi_j(\lambda)\Psi(\lambda).
\]
%Note, however, that this framework is not completely self contained:
%to obtain the evolution of $\p_{Bj}\Theta$ along $\p_{Ai}$ we need to
%know $\p_{Ai+j}\Theta$ 

%=\sum_{i=-j}^{\infty}H_{j\a}\lambda^{\a}
%where the $\lambda$-independent functionals $H_{jk}$ are in involution
%with respect to all Poisson structures ${\cal D}_m$.
 
%The local Poisson structure $\cal D$ can be related to the sequence of 
%symplectic structures ${\bf \Om}^k$ on ${\cal W}_g$; we have that
%Recall that for $\psi,{\psi'}\in {\cal W}_g$ 
%\[
%{\bf \Om}^k(\psi ,\;\psi')=\int_{\delta M} R^k\phi \ast d\psi' -
%\psi' \ast dR^k\psi.
%\]
%and define
%\[
%{\bf \Om}_{\lambda}(\psi ,\;\psi')
%=\sum_{i=0}^{\infty}{\bf \Om}^i(\psi ,\;\psi') \lambda^i ={\bf \Om}((1-\lambda
%R)^{-1}\psi, \psi'). 
%\]
%it is the Poisson structure of 
%${\bf \Om}^0-\lambda {\bf \Om}^1$. Let $\Phi$ be a twistor function.
%From recursion relations
%\[
%({\bf \Om}^0-\lambda {\bf \Om}^1)(\Phi, \phi)=
%{\bf \Om}^0(\Phi, \phi)-\lambda {\bf \Om}^0(R\Phi, \phi)=0
%\]
%so that Hamiltonians defined by the $\Phi$ (treated as a Hamiltonian
%vector fields) are Casimirs for ${\cal D}$.

%Define a current
%\[
%{j_{AA'}}^{(i)}=(0,\;0,\;\p_{\wt}\p _{t_i}\Om,\;\p_{\zt}\p _{t_i}\Om)
%\]
%It follows that
%\be
%\nabla ^{AA'}{j_{AA'}}^{(i)}=\p_{\mu}(g^{\nu \mu}{j_{\nu}}^{(i)})=0
%\ee
%%%%%%%%%%%%%%%%%%%%%%%%%%%%%%%%%%%%%%%%%%%%%%%%%%%%%%%%%%%%%%%%%%%%%%%%%%

\section{Outlook - examples with higher symmetries}
\label{sectFG}
This section motivates the study of solutions to heavenly equations
which are invariant under some hidden symmetries, e.g.\ along the
higher flows.  More generally, one can consider solutions to the
hyper-K\"ahler equations without symmetries, but whose hierarchies do admit
symmetries.

%This should lead to
%a large class of ASD vacuum metrics analogous to the `finite gap'
%solutions in soliton theory \cite{Novikov74}.  
%Imposing three independent hidden symmetries on the twistor space
%should lead to metrics expressible by $\theta$-functions.
%Solutions in terms of
%$\theta$-functions are implicitly given in \cite{DMW95} where ASD
%metrics were expressed in terms of KdV potentials.

In a subsequent paper we shall give a general construction of such
metrics based on a generalisation of \cite{TW79}.  We consider the
case in which the
twistor spaces have a globally defined twistor function
homogeneous of degree $n+1$. This implies that the metric admits a
Killing spinor (some solutions with this property are given by
\cite{D98}).  Global sections $Q\in H^0(\CP^1, \O(n+1))$ on
non-deformed twistor space $\mu :{\cal PT}\longrightarrow \CP^1$ will
be classified and $Q$-preserving deformations of the complex structure
of a neighbourhood of an $\O(1)\oplus\O(1)$ section of $\mu$ will be
studied.  The cohomology classes determining the deformation will
depend on the fibre coordinates of $\mu$ only via $Q$.  The canonical
forms of patching functions can be derived to give explicit solutions
to anti-self-dual ASD vacuum Einstein equation.

There are also further details of the bi-Hamiltonian structure that
could usefully be clarified.

\section {Acknowledgments}
We are grateful to Roger Penrose, George Sparling, Paul Tod, Nick
Woodhouse, and others for some helpful discussions.
Some parts of this work were finished during the workshop
{\em Spaces of geodesics and complex methods in general relativity and 
geometry } held in the summer of 1999 at the Erwin Schr{\"o}dinger 
Institute in Vienna.  We wish to thank ESI  
for the hospitality and for financial
assistance. LJM was supported by NATO grant CRG 950300.

\end{document}